\font\tenmsa=msam10
\font\tenmsb=msbm10

\font\largebf=cmbx10 scaled\magstep2
\def\AGOSS{A^{\cal G}}
\def\all{\hbox{for all}}
\def\and{\hbox{and}}
\def\Abar{\overline{A}}
\def\bra#1#2{\langle#1,#2\rangle}
\def\Bra#1#2{\big\langle#1,#2\big\rangle}
\def\cite#1\endcite{[#1]}
\def\CLB{{\cal CLB}}
\def\dist{\hbox{\rm dist}}
\def\dom{\hbox{\rm dom}}
\def\dbs{^{**}}
\def\episum{\mathop{\nabla}}
\def\eps{\varepsilon}
\def\f#1#2{{#1 \over #2}}
\def\fbar{\overline f}
\def\fourth{\ts\f14}
\def\half{{\textstyle\f12}}

\def\infn{\inf\nolimits}
\def\it{{\wt\iota}}

\def\lr{\Longrightarrow}

\def\minn{\min\nolimits}
\def\on{\hbox{on}}

\def\PC{{\cal PC}}
\def\PCLSC{{\cal PCLSC}}
\def\phi{\varphi}
\def\PQ{{\cal P}_q}
\def\PQT{{\cal P}_{\tilde q}}
\def\pt{\wt p}
\def\qed{\hfill\hbox{\tenmsa\char03}}
\def\qlr{\quad\lr\quad}

\def\qt{\wt q}

\def\quand{\quad\and\quad}
\def\r{\hbox{\tenmsb R}}
\def\rbar{\,]{-}\infty,\infty]}
\def\rl{\Longleftarrow}
\def\rthalf{\ts\f1{\sqrt2}}
\def\rttwo{\sqrt2}
\def\Sem{{\cal S}}
\long\def\slant#1\endslant{{\sl#1}}

\def\supn{\sup\nolimits}
\def\T{{\cal T}}
\def\TCLB{{\cal T_{CLB}}}
\def\TD{{\cal T_D}}
\def\TCLBN{{\cal T_{CLBN}}}
\def\TNORM{{\cal T}_{\|\ \|}}
\def\toto{\ {\mathop{\hbox{\tenmsa\char19}}}\ }
\def\ts{\textstyle}
\def\wh{\widehat}
\def\wt{\widetilde}
\def\xbra#1#2{\lfloor#1,#2\rfloor}
\def\Xbra#1#2{\big\lfloor#1,#2\big\rfloor}
\def\ybra#1#2{\lceil#1,#2\rceil}
\def\Ybra#1#2{\big\lceil#1,#2\big\rceil}
\def\defSection#1{}
\def\defCorollary#1{}
\def\defDefinition#1{}
\def\defExample#1{}
\def\defLemma#1{}
\def\defNotation#1{}
\def\defProblem#1{}
\def\defRemark#1{}
\def\defTheorem#1{}
\def\locno#1{}
\def\meqno#1{\eqno(#1)}
\def\nmbr#1{}
\def\Proof{\medbreak\noindent{\bf Proof.}\enspace}
\def\Proo{{\bf Proof.}\enspace}
\def\Signoff{}
\def \INTsec{1}
\def \SSDsec{2}
\def \PCdef{2.1}
\def \QPOSdef{2.2}
\def \QPOSone{1}
\def \QPOStwo{2}
\def \QPOSthree{3}
\def \QPOSfour{4}
\def \QPOSfive{5}
\def \Hex{2.3}
\def \EEex{2.4}
\def \NOTSSDex{2.5}
\def \ROOTlem{2.6}
\def \ROOTrem{2.7}
\def \POSFdef{2.8}
\def \POSFone{6}
\def \Llem{2.9}
\def \FATdef{2.10}
\def \FATone{7}
\def \PHlem{2.11}
\def \FFATlem{2.12}
\def \FFATrem{2.13}
\def \PHIMAXthm{2.14}
\def \EXTsec{3}
\def \THETAdef{3.1}
\def \LINKone{8}
\def \LINKtwo{9}
\def \LINKthree{10}
\def \LINKfour{11}
\def \LINKfive{12}
\def \LINKsix{13}
\def \THPSlem{3.2}
\def \MAXthm{3.3}
\def \MAXcor{3.4}
\def \MAXone{14}
\def \MAXtwo{15}
\def \MAXCONcor{3.5}
\def \HHATone{16}
\def \HHATtwo{17}
\def \HHATthree{18}
\def \HOMsec{4}
\def \DUALITYdef{4.1}
\def \DUALone{19}
\def \DUALtwo{20}
\def \DUALthree{21}
\def \DUALfour{22}
\def \DUALfive{23}
\def \LDlem{4.2}
\def \HOMlem{4.3}
\def \GOSSEZdef{4.4}
\def \GOSSEZone{24}
\def \AGOSSthm{4.5}
\def \AGOSSone{25}
\def \BSSDsec{5}
\def \NORMdef{5.1}
\def \IOTAtwo{26}
\def \IOTAone{27}
\def \IOTAthree{28}
\def \IOTAfour{29}
\def \IOTAfive{30}
\def \NORMrem{5.2}
\def \EENex{5.3}
\def \VZdef{5.4}
\def \VZone{31}
\def \VZtwo{32}
\def \VZthree{33}
\def \VZfour{34}
\def \Pdef{5.5}
\def \VZEXthm{5.6}
\def \VZEXone{35}
\def \VZEXtwo{36}
\def \EXtwo{37}
\def \EXthree{38}
\def \EXfour{39}
\def \EXfive{40}
\def \EXseven{41}
\def \EXeight{42}
\def \VZEQthm{5.7}
\def \EXthm{5.8}
\def \BSTARrem{5.9}
\def \IOTAnew{43}
\def \Hcor{5.10}
\def \VZONLYthm{5.11}
\def \VZONLYone{44}
\def \VZONLYtwo{45}
\def \VZONLYthree{46}
\def \FPHIrem{5.12}
\def \FPHIone{47}
\def \FPHItwo{48}
\def \TWOrem{5.13}
\def \VZONLYrem{5.14}
\def \DBSSDsec{6}
\def \BSSDdef{6.1}
\def \BSSDone{49}
\def \BSSDthree{50}
\def \BSSDfour{51}
\def \FACTORlem{6.2}
\def \Hrem{6.3}
\def \PSEUDOdef{6.4}
\def \PSEUDOone{52}
\def \EEDUALex{6.5}
\def \EEfive{53}
\def \EEeight{54}
\def \OTHERrem{6.6}
\def \PRODrem{6.7}
\def \RSUMlem{6.8}
\def \SHARPrem{6.9}
\def \SSDDlem{6.10}
\def \Hone{55}
\def \MASdef{6.11}
\def \MASVZthm{6.12}
\def \USEdef{6.13}
\def \USErem{6.14}
\def \EQthm{6.15}
\def \LIMthree{56}
\def \LIMfour{57}
\def \CLBsec{7}
\def \TCBlem{7.1}
\def \PATHOrem{7.2}
\def \TCLBElem{7.3}
\def \QTsec{8}
\def \TDdef{8.1}
\def \TDone{58}
\def \TDtwo{59}
\def \TWOTOPlem{8.2}
\def \EESTARthm{8.3}
\def \BONUSrem{8.4}
\def \MONsec{9}
\def \ABARone{60}
\def \TYPEDdef{9.1}
\def \NIdef{9.2}
\def \NIone{61}
\def \NItwo{62}
\def \TYPEEDdef{9.3}
\def \NIthree{63}
\def \RECASTlem{9.4}
\def \SRone{64}
\def \EEEQthm{9.5}
\def \SRdef{9.6}
\def \SRthm{9.7}
\def \ZAGone{65}
\def \ZAGprob{9.8}
\def \BRthm{9.9}
\def \NAPAIRone{66}
\def \DFPthm{9.10}
\def \FMsec{10}
\def \FMthm{10.1}
\def \FMone{67}
\def \FMtwo{68}
\def \FENMORSUFFone{69}
\def \FENMORSUFFtwo{70}
\def \FENMORSUFFthree{71}
\def \FENMORSUFFfour{72}
\def \BS{1}
\def \FITZ{2}
\def \GOSSEZ{3}
\def \GOSSEZC{4}
\def \ASTWO{5}
\def \ASTHREE{6}
\def \ASD{7}
\def \MT{8}
\def \MLT{9}
\def \PENOT{10}
\def \PRAGUE{11}
\def \FENCHEL{12}
\def \RTRCA{13}
\def \RANGE{14}
\def \MANDM{15}
\def \BR{16}
\def \FIVE{17}
\def \PANDM{18}
\def \HBL{19}
\def \HBM{20}
\def \NRSSD{21}
\def \VZ{22}
\def \ZAGRODNY{23}
\def \ZBOOK{24}
%
\magnification 1200
\parindent0pt
\headline{\ifnum\folio=1
{\hfil{\largebf Banach SSD spaces and classes of monotone sets}\hfil}
\else\centerline{\rm {\bf Banach SSD spaces and classes of monotone sets}}\fi}
\bigskip
\centerline{\bf S. Simons}
\bigskip
\centerline{\sl Department of Mathematics, University of California}
\centerline{\sl Santa Barbara, CA 93106-3080, U.S.A.}
\centerline{\sl simons@math.ucsb.edu}
\medskip
\noindent
In this paper, we unify the theory of SSD spaces and the theory of strongly representable sets, and we apply our results to the theory of the various classes of maximally monotone sets.   In particular, we prove that type (ED), dense type, type (D), type (NI) and strongly representable are equivalent concepts and, consequently, that the known properties of\break strongly representable sets follow from known properties of sets of type (ED).
\defSection \INTsec
\medbreak
\leftline{\bf \INTsec.\quad Introduction}
\medskip
\noindent
In Sections \SSDsec--\HOMsec, we give a more complete version of the algebraic theory of SSD spaces, as introduced in \cite\PANDM\endcite,  and further developed in \cite\HBM\endcite.   Apart from the fact that we write ``$\PQ$'' instead of ``pos'', we use the notation of the latter of these references.   What distinguishes the three sections is the number of bilinear forms considered: one in Section \SSDsec, two in Section \EXTsec, and three in Section \HOMsec.   The concepts of \slant SSD space\endslant, \slant$q$--positive set\endslant\ and the convex function, $\Phi_A$, associated with a $q$--positive set, $A$, are introduced in Definition \QPOSdef.   The functions $\Phi$ are the generalizations to SSD spaces of the Fitzpatrick functions of monotone sets.   In Definition \POSFdef, we introduce the $q$--positive set, $\PQ(f)$, associated with an appropriate convex function, $f$, and in Definition \FATdef, we introduce the \slant intrinsic conjugate\endslant, $f^@$, of a convex function, $f$.   The main result in this section is Theorem \PHIMAXthm, though the results marked ``Lemma'' will be used throughout the paper.   In Section \EXTsec, we consider the situation of an SSD space linked by a linear map to an external vector space.  We then define another convex function, $\Psi_A$, associated with a $q$--positive set, $A$, using a convex function, $\Theta_A$, on the external space.   The basic properties of these functions are collected together in Lemma \THPSlem. In Theorem \MAXthm\ and Corollaries \MAXcor\ and \MAXCONcor, we discuss a maximal property of the $\Psi$ functions, which complements a minimal property of the $\Phi$ functions in certain circumstances.   While the material of Section \SSDsec\ is essentially algebraic \big(apart from the disguised differentiability argument of Lemma \FFATlem(a)\big), Theorem \MAXthm\ uses the Fenchel--Moreau theorem from convex analysis, for a (possibly nonhausdorff) locally convex space.   Since we have not seen this result in the literature, we give a proof of it for the convenience of the reader in Section \FMsec.   In Section \HOMsec, we consider the special situation where the external space is also an SSD space, and the allied concept of \slant SSD--homomorphism\endslant.   These are introduced in Definition \DUALITYdef.  These ideas allow us to apply the analysis of Sections \SSDsec\ and \EXTsec, with the SSD space replaced by the external space.  This enables us to generalize (in Definition \GOSSEZdef\ and Theorem \AGOSSthm) to SSD spaces  some concepts due to Gossez for maximally monotone multifunctions.
\smallskip
In Sections \BSSDsec\ and \DBSSDsec, we specialize to the situation in which the SSD spaces have a Banach space structure also.   In Definition \NORMdef, we introduce the concept of a \slant Banach SSD space\endslant, which is an SSD space with a Banach space structure satisfying the compatibility condition (\IOTAtwo), from which it follows that the Banach space dual can be considered as a linked external space (see Remark \BSTARrem).   Section \BSSDsec\ is inspired by Voisei--Z\u alinescu, \cite\VZ\endcite.    
In Definition \VZdef, we introduce the concept of \slant VZ function\endslant\ on a Banach SSD space.   Our main result on VZ functions, established in Theorems \VZEXthm(d) and \EXthm(b), is that \slant if $f$ is a lower semicontinuous VZ function then $\PQ(f)$ is {\rm maximally} $q$--positive, $f^@$ is also a VZ function, and $\PQ\big(f^@\big) = \PQ(f)$\endslant.   The argument in Theorem \VZEXthm\ uses completeness heavily, as well as the fact that a proper, convex, lower semicontinuous function on a Banach space dominates a continuous affine function.   Theorem \VZEQthm\ contains a characterization of VZ functions in terms of the concept of \slant $p$--density\endslant\ introduced in Definition \Pdef.   We show in Theorem \VZONLYthm\ that if $f$ is a lower semicontinuous VZ function on a Banach SSD space and $A = \PQ(f)$ then $f$ lies between $\Psi_A$ and $\Phi_A$ and, further, if $h$ is a proper convex function on $B$ that lies between $\Psi_A$ and $\Phi_A$ then $h$ and $h ^@$ are VZ functions.   The rest of Section \BSSDsec\ is devoted to some counterexamples.   Section \BSSDsec\ uses the material of Sections \SSDsec\ and \EXTsec, but not Section \HOMsec.   By contrast, Section \DBSSDsec\ depends heavily on Section \HOMsec.  Here we consider the situation where the dual of a Banach SSD space can be made into a Banach SSD space in its own right, satisfying the compatibility condition (\BSSDone), and we write $\qt$ for the function on the dual corresponding to the function $q$ previously defined on the original Banach space.   In Definition \MASdef, we work towards defining an analog for SSD spaces of the concept of \slant strongly representable multifunction\endslant, as expounded by Voisei and Z\u alinescu in \cite\VZ\endcite\ and Marques Alves and Svaiter in \cite\ASTWO\endcite, \cite\ASTHREE\endcite\ and \cite\ASD\endcite.   In order to to this, we introduce the concept of \slant MAS function\endslant\ in Definition \MASdef.   The first main result of Section \DBSSDsec\ is Theorem \MASVZthm\ (which leads to Theorem \SRthm), in which we establish that, under a certain mild side condition, the concepts of MAS function and VZ function are identical.   The main tools here are Rockafellar's formula for the conjugate of the sum of convex functions, and the fact that the conjugate of the function $\half\|\cdot\|^2$ on a Banach space is the function $\half\|\cdot\|^2$ on the dual space, which are both used in Lemma \SSDDlem.   The other main result of Section \DBSSDsec\ is Theorem \EQthm\ (which leads to Theorem \EEEQthm), which depends on the concept of a \slant compatible\endslant\ topology on $B^*$ introduced in Definition \USEdef, and describes the relationship between compatible topologies and the Gossez extension of a maximally $q$--positive set introduced in Definition \GOSSEZdef.    
\smallskip
Sections \CLBsec\ and \QTsec\ are devoted to a discussion of certain esoteric topologies on the bidual of a Banach space and the Banach SSD dual of a Banach SSD space.   Here is some background for the problem.   Suppose that $E$ is a nonzero Banach space, and consider the function $\qt\colon\ (x^*,x\dbs) \mapsto \bra{x^*}{x\dbs}$ from $(E \times E^*)^* = E^* \times E\dbs$ into $\r$.   While it is true that the norm topology on $E^* \times E\dbs$ makes $\qt$ continuous, it has been known since the work of Gossez on maximally monotone multifunctions that the norm topology is too large to be of any practical use.   The reason for this can be traced to the fact that it is not generally compatible in the sense of Definition \USEdef.   (See Remark \USErem.)  In Section \QTsec, we introduce the topologies $\TD$ on Banach SSD duals, which have the properties that they are sufficiently small that they are compatible and sufficiently large that Theorem \EQthm\ leads to significant results.
The topologies $\TD$ are based on the topologies $\TCLB$ on the bidual of a Banach space that have been previously studied, the properties of which are stated in Section \CLBsec.
\smallskip
So far, we have been describing general theories, but we have not discussed any particular examples.   In Example \Hex, we give three examples of SSD spaces, of which (c) is probably the most interesting.  As we observe in Remark \NORMrem, Example \Hex(a,b,c) are actually Banach SSD spaces.   Example \EEex\ is the example that leads to results on monotonicity --- it is shown in Example \EENex\ how to norm this example so that it becomes a Banach SSD space, and in Example \EEDUALex\ how to make its dual a Banach SSD dual.   We then show in Theorem \EESTARthm\ that $\qt$ is continuous with respect to the topology $\TD$, and so we are now in the position that we can apply Theorems \VZEXthm, \MASVZthm\ and \EQthm\ to this example.   This leads to the results on monotonicity that we give in Section \MONsec.
\smallskip
We start Section \MONsec\ with a brief history of various classes of maximally monotone sets.   We define \slant type (D)\endslant, \slant dense type\endslant, \slant type (NI)\endslant, \slant type (WD)\endslant, \slant type (ED)\endslant\ and \slant strongly representable\endslant\ in Definitions \TYPEDdef, \NIdef, \TYPEEDdef\ and \SRdef.   The easy implications are that, for maximally monotone sets,\quad type (ED) $\lr$ dense type $\lr$ type (D) $\lr$ type (WD) $\lr$ type (NI).\quad   Marques Alves and Svaiter proved recently in \cite\ASD\endcite\ that\quad type (NI) $\lr$ type (D).\quad In Theorem \EEEQthm, we show how the techniques discussed in this paper lead to the stronger conclusion that\quad type (NI) $\lr$ type (ED).\quad   The obvious significance of this is that Theorem \EEEQthm\ leads to solutions of several problems that have been open for some time.   These issues are discussed in the paragraph preceding Theorem \EEEQthm.   However, Theorem \EEEQthm\ is significant for another reason.   Strongly representable sets were initially introduced in \cite\ASTWO\endcite\ and \cite\VZ\endcite, and it was proved in \cite\ASTHREE\endcite\ that a set is strongly representable $\iff$ it is maximally monotone of type (NI).   In Theorem \SRthm, we show how the techniques discussed in this paper lead to a proof of this equivalence.   If we now combine together the results discussed above, we see that a set is strongly representable $\iff$ it is maximally monotone of type (ED).   This enables us to use the properties known for maximally monotone sets of type (ED) to obtain results about strongly representable sets, which frequently improve on the results already known.   We give these results in Theorems \BRthm\ and \DFPthm, with references to what is currently in the literature.                 
\smallskip
In the Appendix, Section \FMsec, we give a proof of the  Fenchel--Moreau theorem for a (possibly nonhausdorff) locally convex space, which we used in Theorem \MAXthm.
\smallskip
The author would like to thank Constantin Z\u alinescu for making him aware of the preprints \cite\ASTWO\endcite\ and \cite\VZ\endcite, and Maicon Marques Alves and Benar Svaiter for making him aware of the preprints \cite\ASTHREE\endcite\ and \cite\ASD\endcite.   He would also like to thank Radu Ioan Bo\c t and Constantin Z\u alinescu for some very perceptive comments on earlier versions of this paper.   The author has learned that the acronym ``SSD'' has also been used to stand for ``strongly subdifferentiable''.   He hopes that this will not cause any confusion.  Finally, he would like to thank the anonymous referees, whose insightful comments resulted in great improvements. 
\defSection\SSDsec
\medbreak
\leftline{\bf \SSDsec.\quad SSD spaces}
\defDefinition \PCdef
\medbreak
\noindent
{\bf Definition \PCdef.}\enspace If $X$ is a nonzero real vector space and $f\colon\ X \to \rbar$, we write $\dom\,f$ for the set $\big\{x \in X\colon\ f(x) \in \r\big\}$.   $\dom\,f$ is the \slant effective domain \endslant of $f$.   We say that $f$ is \slant proper\endslant\ if $\dom\,f \ne \emptyset$.   We write $\PC(X)$ for the set of all proper convex functions from $X$ into $\rbar$.   If $X$ is a nonzero real Banach space, we write $X^*$ for the dual space of $X$ \big(with the pairing $\bra\cdot\cdot\colon X \times X^* \to \r$\big).   
\defDefinition \QPOSdef
\medbreak
\noindent
{\bf Definition \QPOSdef.}\enspace We will say that $\big(B,\xbra\cdot\cdot\big)$ is a \slant symmetrically self--dual space (SSD space)\endslant\  if $B$ is a nonzero real vector space and $\xbra\cdot\cdot\colon B \times B \to \r$ is a symmetric bilinear form.   In this case, we will \slant always\endslant\ write\quad $q(b) := \half\xbra{b}{b}$\quad($b \in B$).   (``$q$'' stands for ``quadratic''.)\quad   We do not assume that $\xbra\cdot\cdot$ separates the points of $B$, as was done in  \cite\PANDM\endcite\ and \cite\HBM\endcite.  With this caveat, the definitions and many of the results of his section appear in \cite\PANDM\endcite\ and \cite\HBM\endcite. 
\smallskip
Now let $\big(B,\xbra\cdot\cdot\big)$ be an SSD space and $A \subset B$.   We say that $A$ is \slant$q$--positive\endslant\ if $A \ne \emptyset$ and
$$b,c \in A \lr q(b - c) \ge 0.$$
In this case, since $q(0) = 0$,
$$b \in A \lr \inf q(A - b) = 0.\meqno\QPOSone$$
We then define $\Phi_A\colon\ B \to \rbar$ by
$$\Phi_A(b) := \supn_A\big[\xbra\cdot{b} - q\big]\quad(b \in B).\meqno\QPOStwo$$
$\Phi_A$ is a generalization to SSD spaces of the ``Fiztpatrick function'' of a monotone set, which was originally introduced in \cite\FITZ\endcite\ in 1988, but lay dormant until it was rediscovered by  Mart\'\i nez-Legaz and Th\'era in \cite\MLT\endcite\ in 2001.
We note then that, for all $b \in B$,
$$\left.\eqalign{
\Phi_{A}(b) &= q(b) - \infn_{a \in A}\big[q(a) - \xbra{a}{b} + q(b)\big]\cr
&= q(b) - \infn_{a \in A}q(a - b) = q(b) - \inf q(A - b).}\right\}\meqno\QPOSthree$$
From (\QPOSone),
$$\Phi_A = q\ \on\ A.\meqno\QPOSfour$$
Thus $\Phi_A \in \PC(B)$.   We say that $A$ is \slant maximally $q$--positive\endslant\ if $A$ is $q$--positive and $A$ is not properly contained in any other $q$--positive set.   In this case, if $b \in B$ and $\inf q(A - b) \ge 0$ then clearly $b \in A$.   In other words,\qquad \big($b \in B \setminus A \lr \inf q(A - b) < 0$\big).\qquad From (\QPOSone),\qquad $\inf q(A - b) \le 0$\qquad and\qquad  \big($\inf q(A - b) = 0 \iff b \in A$\big).\qquad Thus, from (\QPOSthree)
$$\Phi_A \ge q\ \on\ B\quand\big(\Phi_A(b) = q(b) \iff b \in A\big).\meqno\QPOSfive$$
We make the elementary observation that if $b \in B$ and $q(b) \ge 0$ then the linear span $\r b$ of $\{b\}$ is $q$--positive.
\medbreak
We now give some examples of SSD spaces and their associated $q$--positive sets.   These examples are taken from \cite\HBM, pp.\ 79--80\endcite.   
\defExample \Hex 
\medbreak
\noindent
{\bf Example \Hex.}\enspace Let $B$ be a Hilbert space with inner product $(b,c) \mapsto \bra{b}{c}$ and $T\colon B \to B$ be a self--adjoint linear operator.   Then $\big(B,\xbra\cdot\cdot\big)$ is an SSD space with $\xbra{b}{c} := \bra{b}{Tc}$ and $q(b) = \half\bra{Tb}{b}$.   Here are three special cases of this example:
\smallbreak
(a)\enspace If, for all $b \in B$, $Tb = b$ then $\xbra{b}{c} := \bra{b}{c}$, $q(b) = \half\|b\|^2$ and every subset of $B$ is $q$--positive
\smallbreak
(b)\enspace If, for all $b \in B$, $Tb = -b$ then $\xbra{b}{c} := -\bra{b}{c}$, $q(b) = -\half\|b\|^2$ and the $q$--positive sets are the singletons.
\smallbreak
(c)\enspace If $B = \r^3$ and $T(b_1,b_2,b_3) = (b_2,b_1,b_3)$ then
$$\Xbra{(b_1,b_2,b_3)}{(c_1,c_2,c_3)} := b_1c_2 + b_2c_1 + b_3c_3,$$
and\quad$q(b_1,b_2,b_3) = b_1b_2 + \half b_3^2$.   Here, if $M$ is any nonempty monotone subset of $\r \times \r$ (in the obvious sense) then $M \times \r$ is a $q$--positive subset of $B$.   The set $\r(1,-1,2)$ is a $q$--positive subset of $B$ which is not contained in a set $M \times \r$ for any monotone subset of $\r \times \r$.   The helix $\big\{(\cos\theta,\sin\theta,\theta)\colon \theta \in \r\big\}$ is a $q$--positive subset of $B$, but if $0 < \lambda < 1$ then the helix $\big\{(\cos\theta,\sin\theta,\lambda\theta)\colon \theta \in \r\big\}$ is not.
\defExample \EEex 
\medbreak
\noindent
{\bf Example \EEex.}\enspace Let $E$ be a nonzero Banach space and $B := E \times E^*$.   For all $(x,x^*)$ and $(y,y^*) \in B$, we set $\Xbra{(x,x^*)}{(y,y^*)} := \bra{x}{y^*} + \bra{y}{x^*}$.  Then $\big(B,\xbra\cdot\cdot\big)$ is an SSD space with $q(x,x^*) = \half\big[\bra{x}{x^*} + \bra{x}{x^*}\big] = \bra{x}{x^*}$.   Consequently, if $(x,x^*), (y,y^*) \in B$ then
$\bra{x - y}{x^* - y^*} = q(x - y,x^* - y^*) = q\big((x,x^*) - (y,y^*)\big)$.   
Thus if $A \subset B$ then $A$ is $q$--positive exactly when $A$ is a nonempty monotone subset of $B$ in the usual sense, and $A$ is maximally $q$--positive exactly when $A$ is a maximally monotone subset of $B$ in the usual sense.   We point out that any finite dimensional SSD space of the form described here must have \slant even\endslant\ dimension.  Thus cases of Example \Hex\ with finite odd dimension cannot be of this special form.
\defExample \NOTSSDex
\medbreak
\noindent
{\bf Example \NOTSSDex.}\enspace $\big(\r^3,\xbra\cdot\cdot\big)$ is \slant not\endslant\ an SSD space with
$$\Xbra{(b_1,b_2,b_3)}{(c_1,c_2,c_3)} := b_1c_2 + b_2c_3 + b_3c_1.$$
(The bilinear form $\xbra\cdot\cdot$ is not symmetric.)
\defLemma \ROOTlem
\medbreak
\noindent
{\bf Lemma \ROOTlem.}\enspace\slant Let $\big(B,\xbra\cdot\cdot\big)$ be an SSD space, $f \in \PC(B)$, $f \ge q$ on $B$ and $b,c \in B$.   Then
$$-q(b - c) \le \Big[\sqrt{(f - q)(b)} + \sqrt{(f - q)(c)}\Big]^2.$$\endslant
\Proof We can and will suppose that $0 \le (f - q)(b)  < \infty$ and $0 \le (f - q)(c)  < \infty$.   Let $\sqrt{(f - q)(b)} < \beta < \infty$ and $\sqrt{(f - q)(c)} < \gamma < \infty$, so that $\beta^2 + q(b) > f(b)$ and $\gamma^2 + q(c) > f(c)$.   Then, writing $\alpha := \beta + \gamma$,
$$\eqalign{\beta\gamma + (\gamma q(b) + \beta q(c))/\alpha
&= \gamma\big(\beta^2 + q(b)\big)/\alpha + \beta\big(\gamma^2 + q(c)\big)/\alpha\cr
&> \gamma f(b)/\alpha + \beta f(c)/\alpha \ge
f\big({(\gamma b + \beta c)}/\alpha\big)\cr
&\ge q\big({(\gamma b + \beta c)}/\alpha\big)
= \big(\gamma^2q(b) + \gamma\beta\xbra{b}{c} + \beta^2q(c)\big)/\alpha^2.}$$
Clearing of fractions, we obtain
$$\alpha^2\beta\gamma + \alpha\big(\gamma q(b) + \beta q(c)\big) >
\gamma^2q(b) + \gamma\beta\xbra{b}{c} + \beta^2q(c),$$
from which $\alpha^2\beta\gamma > -\beta\gamma q(b) + \beta\gamma\xbra{b}{c} - \beta\gamma q(c) = -\beta\gamma q(b - c)$.   If we now divide by $\beta\gamma$, we obtain $\alpha^2 > -q(b - c)$, and the result follows by letting $\beta \to \sqrt{(f - q)(b)}$ and $\gamma \to \sqrt{(f - q)(c)}$.\qed 
\defRemark \ROOTrem
\medbreak
\noindent
{\bf Remark \ROOTrem.}\enspace It follows from Lemma \ROOTlem\ and the Cauchy--Schwarz inequality that
$$-q(b - c) \le 2(f - q)(b) + 2(f - q)(c).$$
In the situation of Example \EEex, we recover \cite\VZ, Proposition 1\endcite.
\defDefinition \POSFdef
\medbreak
\noindent
{\bf Definition \POSFdef.}\enspace If $\big(B,\xbra\cdot\cdot\big)$ is an SSD space, $f \in \PC(B)$ and $f \ge q$ on $B$, we write
$$\PQ(f) := \big\{b \in B\colon\ f(b) = q(b)\big\}.$$
We note then that  (\QPOSfive) implies that
$$\hbox{\sl if}\ A\ \hbox{\sl is maximally}\ q\hbox{\sl--positive then}\ A = \PQ(\Phi_A).\meqno\POSFone$$
\par
The following result is suggested by Burachik--Svaiter, \cite\BS, Theorem 3.1, pp. 2381--2382\endcite\ and Penot, \cite\PENOT, Proposition 4(h)$\lr$(a), pp. 860--861\endcite.
\defLemma \Llem
\medbreak
\noindent
{\bf Lemma \Llem.}\enspace\slant Let $\big(B,\xbra\cdot\cdot\big)$ be an SSD space, $f \in \PC(B)$, $f \ge q$ on $B$ and $\PQ(f) \ne \emptyset$.   Then $\PQ(f)$ is a $q$--positive subset of $B$.\endslant
\Proof This is immediate from Lemma \ROOTlem.\qed
\medbreak
We now introduce a concept of conjugate that is intrinsic to an SSD space without any topological conditions.
\defDefinition \FATdef
\medbreak
\noindent
{\bf Definition \FATdef.}\enspace If $\big(B,\xbra\cdot\cdot\big)$ is an SSD space and $f \in \PC(B)$, we write $f^@$ for the Fenchel conjugate of $f$ with respect to the pairing $\xbra\cdot\cdot$, that is to say,
$$\all\ c \in B,\qquad f^@(c) := \supn_B\big[\xbra{\cdot}{c} - f\big].\meqno\FATone$$
\par
The next result gives some basic properties of ${\Phi_A}^@$ and ${\Phi_A}^{@@}$.   It will be used in Theorem \PHIMAXthm, Lemma \LDlem(e) and Lemma \HOMlem(c).   The proof of Lemma \PHlem(c) below is due to Radu Ioan Bo\c t.
\defLemma \PHlem
\medbreak
\noindent
{\bf Lemma \PHlem.}\enspace\slant Let $\big(B,\xbra\cdot\cdot\big)$ be an SSD space and $A$ be a nonempty $q$--positive subset of $B$.   Then:
\smallbreak
\noindent
{\rm(a)}\enspace ${\Phi_A}^@ \le q$ on $A$.
\smallbreak
\noindent
{\rm(b)}\enspace ${\Phi_A}^@ \ge  \Phi_A \vee q$ on $B$.
\smallbreak
\noindent
{\rm(c)}\enspace ${\Phi_A}^{@@} = \Phi_A$ on $B$.
\endslant
\Proof Let $a \in A$ and $b \in B$.   From (\QPOStwo), $\xbra{a}{b} - q(a) \le \Phi_A(b)$.   Thus $\xbra{b}{a} - \Phi_A(b) \le q(a)$, and we obtain (a) by taking the supremum over $b \in B$.   Let $c \in B$.   Then, from (\QPOSfour), 
$$\eqalignno{{\Phi_A}^@(c) 
&= \supn_B\big[\xbra{c}\cdot - \Phi_A\big]
\ge \big[\xbra{c}{c} - \Phi_A(c)\big] \vee \supn_A\big[\xbra{c}\cdot - \Phi_A\big]\cr
&= \big[2q(c) - \Phi_A(c)\big] \vee \supn_A\big[\xbra{c}\cdot - q\big]
= \big[2q(c) - \Phi_A(c)\big] \vee \Phi_A(c).}$$
Now if $\Phi_A(c) = \infty$ then obviously $\big[2q(c) - \Phi_A(c)\big] \vee \Phi_A(c) \ge q(c)$, while if $\Phi_A(c) \in \r$ then $\big[2q(c) - \Phi_A(c)\big] \vee \Phi_A(c) \ge \half[2q(c) - \Phi_A(c)\big] + \half\Phi_A(c) = q(c)$.   Thus ${\Phi_A}^@(c) \ge  \Phi(c) \vee q(c)$. This completes the proof of (b).   From (a), for all $b \in B$, ${\Phi_A}^{@@}(b) = \supn_B\big[\xbra\cdot{b} - {\Phi_A}^@\big]\ge \supn_A\big[\xbra{b}\cdot - {\Phi_A}^@\big] \ge \supn_A\big[\xbra{b}\cdot - q\big] =\Phi_A(b)$.   Thus\quad ${\Phi_A}^{@@} \ge \Phi_A$ on $B$.\quad However, it is obvious that\quad ${\Phi_A}^{@@} \le \Phi_A$ on $B$,\quad which completes the proof of (c).\qed
\medbreak
Our next result represents an improvement of the result proved in \cite\HBM, Lemma 19.12, p.\ 82\endcite, and uses a disguised differentiability argument.   Lemma \FFATlem\ will be used in Theorem \PHIMAXthm, Corollary \MAXcor, Lemma \LDlem\ and Theorem \EXthm(b). See Remark \FFATrem\ below for another proof of Lemma \FFATlem(a), due to Constantin Z\u alinescu.   
\defLemma \FFATlem
\medbreak
\noindent
{\bf Lemma \FFATlem.}\enspace\slant Let $\big(B,\xbra\cdot\cdot\big)$ be an SSD space, $f \in \PC(B)$ and  $f \ge q$ on $B$.   Then:
$$\leqalignno{a \in \PQ(f)\ \and\ b \in B &\qlr \xbra{b}{a} \le q(a) + f(b).&{\rm(a)}\cr
f^@ &= q\ \on\ \PQ(f).&{\rm(b)}\cr
\hbox{If}\ \PQ(f) \ne \emptyset\quad &\hbox{then}\quad f \ge \Phi_{\PQ(f)}\ \on\ B.&{\rm(c)}}$$
\endslant
\Proo Let $a \in \PQ(f)$ and $b \in B$.   Let $\lambda \in \,]0,1[\,$.   For simplicity in writing, let $\mu := 1 - \lambda \in \,]0,1[\,$.   Then
$$\eqalign{\lambda^2q(b) + \lambda\mu\xbra{b}{a} + \mu^2q(a) &= q\big(\lambda b + \mu a\big) \le f(\lambda b + \mu a)\cr
&\le \lambda f(b) + \mu f(a) = \lambda f(b) + \mu q(a).}$$
Thus\quad $\lambda^2q(b) + \lambda\mu\xbra{b}{a} \le \lambda f(b) + \lambda\mu q(a)$.\quad We now obtain (a) by dividing by $\lambda$ and letting $\lambda \to 0$.   Now let $a \in \PQ(f)$.   From (a), \quad$b \in B \lr \xbra{a}{b} - f(b) \le q(a)$,\quad and it follows by taking the supremum over $b \in B$ that $f^@(a) \le q(a)$.   On the other hand,\quad $f^@(a) \ge \xbra{a}{a} - f(a) = 2q(a) - q(a) = q(a)$,\quad completing the proof of (b).   Finally, let $b \in B$ and $a \in \PQ(f)$.   Then, from (a), $f(b)\ge \xbra{b}{a} - q(a)$.   Taking the supremum over $a \in \PQ(f)$ and using (\QPOStwo), $f(b) \ge \Phi_{\PQ(f)}(b)$.   Thus\quad $f \ge \Phi_{\PQ(f)}$ on $B$,\quad giving (c).\qed
\defRemark \FFATrem
\medbreak
\noindent
{\bf Remark \FFATrem.}\enspace The author is grateful to Constantin Z\u alinescu for pointing out to him the following alternative proof of Lemma \FFATlem(a).   From Lemma \ROOTlem, with $c$ replaced by $a$,
$-q(b) + \xbra{b}{a} - q(a) = -q(b - a) \le (f - q)(b)$.   Thus $\xbra{b}{a} - q(a) \le f(b)$, as required. 
\defTheorem  \PHIMAXthm
\medbreak
\noindent
{\bf Theorem  \PHIMAXthm.}\enspace\slant Let $\big(B,\xbra\cdot\cdot\big)$ be an SSD space and $A$ be a maximally $q$--positive subset of $B$.   Then\quad ${\Phi_A}^@ \ge  \Phi_A \ge q$ on $B$\quand $\PQ\big({\Phi_A}^@\big) = \PQ\big(\Phi_A\big) = A$.\endslant
\Proof The first assertion follows from Lemma \PHlem(b) and (\QPOSfive).   It is clear from this and (\POSFone) that $\PQ\big({\Phi_A}^@\big) \subset \PQ\big(\Phi_A\big) = A$.   On the other hand, we can apply Lemma \FFATlem(b) to $f := \Phi_A$ and obtain $\PQ(\Phi_A) \subset \PQ({\Phi_A}^@)$, which gives the second assertion.\qed   
\defSection\EXTsec
\medbreak
\leftline{\bf \EXTsec.\quad SSD spaces with a linked external space}
\medskip
\noindent
A word is in order about the \slant conjugate\endslant\ of a convex function.   If a vector space is paired with itself by a bilinear form, we use the notation $^@$ to denote the conjugate with respect to this pairing.   We have already seen an example of this in (\FATone), and we will see another one in (\DUALfive).   If a vector space $X$ is paired with another vector space $Y$ by a bilinear form $\bra\cdot\cdot\colon\ X \times Y \to \r$ and $f \in \PC(X)$, we write $f^*$ for the conjugate of $f$ with respect to this pairing, that is to say, for all $y \in Y$, $f^*(y) := \sup_X\big[\bra{\cdot}{y} - f\big]$.   We will have an example of this in (\LINKtwo): we will come to another case in the first paragraph of Section \CLBsec.   
\medbreak
We now introduce an important situation, in which an SSD space $\big(B,\xbra\cdot\cdot\big)$ supports a second duality other than that defined by $\xbra\cdot\cdot$.
\defDefinition \THETAdef
\medbreak
\noindent
{\bf Definition \THETAdef.}\enspace Let $\big(B,\xbra\cdot\cdot\big)$ be an SSD space.   We say that $\big(D, \iota,\bra\cdot\cdot\big)$ is a \slant linked external space\endslant\ if $D$ is a nonzero real vector space, $\iota\colon\ B \to D$ is a linear map and\break $\bra\cdot\cdot\colon\ B \times D \to \r$ is a bilinear form such that
$$\all\ b,c\in B,\quad \Bra{b}{\iota(c)} = \xbra{b}{c}.\meqno\LINKone$$
If $\big(B,\xbra\cdot\cdot\big)$ is an SSD space and $\iota$ is the identity map on $B$ then $\big(B, \iota,\xbra\cdot\cdot\big)$ is a linked external space.   We will discuss more interesting examples of totally differing characters in Definition \DUALITYdef\ and Remark \BSTARrem.
\smallskip
Let $\big(D, \iota,\bra\cdot\cdot\big)$ be a linked external space.   If $f \in \PC(B)$ and $d \in D$ then we have
$$f^*(d) := \supn_B\big[\bra\cdot{d} - f\big].\meqno\LINKtwo$$
It is clear from (\LINKtwo), (\LINKone) and (\FATone) that 
if $f \in \PC(B)$ then
$$f^* \circ \iota = f^@.\meqno\LINKthree$$
\par
Now let $A$ be a nonempty $q$--positive subset of $B$.   We then define the function $\Theta_A\colon\ D \to \rbar$ by,
$$\all\ d \in D,\quad \Theta_A(d) := \supn_A\big[\bra\cdot{d} - q\big] = \supn_A\big[\bra\cdot{d} - \Phi_A\big]\meqno\LINKfour$$
\big(the equality of the two expressions follows from (\QPOSfour)\big).   It is clear from the first expresson in (\LINKfour), (\LINKone) and (\QPOStwo) that
$$\Theta_A \circ \iota = \Phi_A,\meqno\LINKfive$$
and so (\QPOSfour) implies that $\Theta_A \in \PC(D)$.   We define the function $\Psi_A\colon\ B \to \rbar$ by
$$\Psi_A := \supn_{d \in D}\big[\bra\cdot{d} - \Theta_A(d)\big].
\meqno\LINKsix$$
\par
In the following lemma, we collect together the basic properties of the functions $\Theta_A$ and $\Psi_A$.   The proof of Lemma \THPSlem(c) below is due to Radu Ioan Bo\c t.   

\defLemma \THPSlem
\medbreak
\noindent
{\bf Lemma \THPSlem.}\enspace\slant Let $\big(B,\xbra\cdot\cdot\big)$ be an SSD space, $\big(D, \iota,\bra\cdot\cdot\big)$ be a linked external space and $A$ be a nonempty $q$--positive subset of $B$.   Then:
\smallskip
\noindent
{\rm(a)}\enspace $\Psi_A \le q$ on $A$. {\rm (Compare Lemma \PHlem(a).)} 
\smallskip
\noindent
{\rm(b)}\enspace $\Psi_A \in \PC(B)$.
\smallskip
\noindent
{\rm(c)}\enspace ${\Phi_A}^* \ge \Theta_A = {\Psi_A}^*$ on $D$. {\rm (Compare Lemma \PHlem(c).)}
\smallskip
\noindent
{\rm(d)}\enspace $\Phi_A = {\Psi_A}^@$.
\smallskip
\noindent
{\rm(e)}\enspace $\Psi_A \ge {\Phi_A}^@ \ge q$ on $B$.
\smallskip\noindent
{\rm(f)}\enspace $A \subset \PQ(\Psi_A) \subset \PQ\big({\Phi_A}^@\big)$.\endslant
\Proof(a)\enspace Let $d \in D$.   Then the first expresson in (\LINKfour) implies that\quad $\bra\cdot{d} - q \le \Theta_A(d)$ on $A$\quad and so\quad $\bra\cdot{d} - \Theta_A(d) \le q$ on $A$.\quad (a) follows by taking the supremum over $d \in D$ and using (\LINKsix), and (b) is immediate from (a).
\smallskip
(c)\enspace We note from (\LINKtwo) and the second expression in (\LINKfour) that, for all $d \in D$,
$${\Phi_A}^*(d) = \supn_B\big[\bra\cdot{d} - \Phi_A\big] \ge \supn_A\big[\bra\cdot{d} - \Phi_A\big] = \Theta_A(d),$$
thus\quad${\Phi_A}^* \ge \Theta_A$ on $D$.\quad From (a) and the second expression in (\LINKfive), for all $d \in D$,
$${\Psi_A}^*(d) = \supn_B\big[\bra\cdot{d} - \Psi_A\big]\ge \supn_A\big[\bra\cdot{d} - \Psi_A\big] \ge \supn_A\big[\bra\cdot{d} - q\big] =\Theta_A(d).$$   Thus\quad ${\Psi_A}^* \ge \Theta_A$ on $B$.\quad However, it is obvious that ${\Psi_A}^* \le \Theta_A$ on $B$,\quad which completes the proof of (c).   (d) follows by composing the equality in (c) with $\iota$ and using (\LINKfive) and (\LINKthree).   (d) implies that ${\Psi_A}^{@@} = {\Phi_A}^@$, and since\quad $\Psi_A \ge {\Psi_A}^{@@}$ on $B$,\quad we obtain the first inequality in (e).   The second inequality in (e) follows from Lemma \PHlem(b). 
\smallskip
(f) is immediate from (a) and (e).\qed
\medbreak
Our next result, which gives a useful maximal property of $\Psi_A$, is motivated by results originally proved by Burachik and Svaiter in \cite\BS\endcite\ for maximally monotone multifunctions.   We will return to this topic in Theorem \VZONLYthm.   We note that the argument used in Theorem \MAXthm\ is similar to that already used in Lemma \THPSlem(c).   
\defTheorem \MAXthm
\medbreak
\noindent
{\bf Theorem \MAXthm.}\enspace\slant Let $\big(B,\xbra\cdot\cdot\big)$ be an SSD space, $\big(D,\iota,\bra\cdot\cdot\big)$ be a linked external space, $A$ be a nonempty $q$--positive subset of $B$, $f \in \PC(B)$ be $w(B,D)$--lower semicontinuous and\quad $f \le q$ on $A$.\quad Then\quad $f^* \ge \Theta_A$ on D\quad and $\Psi_A$ is the largest $w(B,D)$--lower semicontinuous element of $\PC(B)$ that is dominated by $q$ on $A$.\endslant
\Proof From (\LINKtwo) and the first expression in (\LINKfour), for all $d \in D$,
$$f^*(d) =\supn_B\big[\bra\cdot{d} - f\big] \ge \supn_A\big[\bra\cdot{d} - f\big] \ge \supn_A\big[\bra\cdot{d} - q\big] = \Theta_A(d),$$
thus \quad $f^* \ge \Theta_A$ on D,\quad as required.  Since $f$ is $w(B,D)$--lower semicontinuous, \hbox{using} the Fenchel--Moreau theorem for the (possibly nonhausdorff) locally convex space $\big(B,w(B,D)\big)$ (see Theorem \FMthm), what we have already proved, and (\LINKsix),
$$f = \supn_{d \in D}\big[\bra{\cdot}{d} - f^*(d)\big] \le \supn_{d \in D}\big[\bra{\cdot}{d} - \Theta_A(d)\big] = \Psi_A\ \on\ B,$$
thus\quad $f \le \Psi_A$ on $B$.\quad  The result follows from Lemma \THPSlem(a), since $\Psi_A$ is obviously $w(B,D)$--lower semicontinuous.\qed     
\defCorollary \MAXcor
\medbreak
\noindent
{\bf Corollary \MAXcor.}\enspace\slant Let $\big(B,\xbra\cdot\cdot\big)$ be an SSD space, $\big(D,\iota,\bra\cdot\cdot\big)$ be a linked external space, $f \in \PC(B)$ be $w(B,D)$--lower semicontinuous\quad $f \ge q$ on $B$,\quad and $A := \PQ(f) \ne \emptyset$.   Then
$$\Psi_A \ge f \ge \Phi_A\ \on\ B \quand {\Phi_A}^* \ge f^* \ge \Theta_A\ \on\ D.\meqno\MAXone$$
If, further, $A$ is maximally $q$--positive then
$$\PQ\big(\Psi_A\big) = \PQ\big({\Phi_A}^@\big) = \PQ\big(\Phi_A\big) = A.\meqno\MAXtwo$$\endslant
\Proof It is clear from Lemmas \Llem\ and \FFATlem(c) that $A$ is a $q$--positive subset of $B$ and\quad $f \ge \Phi_A$ on $B$,\quad from which\quad ${\Phi_A}^* \ge f^*$ on $D$.\quad Theorem \MAXthm\ implies that\quad $f^* \ge \Theta_A$ on $D$\quad and\quad $\Psi_A \ge f$ on $B$,\quad which completes the proof of (\MAXone).   If we use Theorem \PHIMAXthm\ to strengthen Lemma \THPSlem(f), we obtain $A \subset \PQ(\Psi_A) \subset \PQ({\Phi_A}^@) = \PQ(\Phi_A) = A$, which completes the proof of (\MAXtwo).\qed
\medbreak
Our next result is a partial converse to Corollary \MAXcor.   We note that $f$ is \slant not\endslant\ required to be $w(B,D)$--lower semicontinuous.  
\defCorollary \MAXCONcor
\medbreak
\noindent
{\bf Corollary \MAXCONcor.}\enspace\slant Let $\big(B,\xbra\cdot\cdot\big)$ be an SSD space, $\big(D, \iota,\bra\cdot\cdot\big)$ be a linked external space, $A$ be a maximally $q$--positive subset of $B$, $f \in \PC(B)$ and\quad
$$\Psi_A \ge f \ge \Phi_A\ \on\ B.\meqno\HHATone$$
Then\quad $f \ge q$ on $B$,\quad $f^@ \ge q$ on $B$\quad and
$$\PQ(f) = \PQ\big(f^@\big) = A = \PQ\big(\Psi_A\big) = \PQ\big({\Phi_A}^@\big) = \PQ\big(\Phi_A\big).$$
\endslant
\Proo From (\HHATone) and (\QPOSfive),\quad$\Psi_A \ge f \ge \Phi_A \ge q\ \on\ B$.\quad Thus, using Lemma \THPSlem(f),
$$A \subset \PQ(\Psi_A) \subset \PQ(f) \subset \PQ(\Phi_A).\meqno\HHATtwo$$
Taking conjugates in (\HHATone) and using Lemma \THPSlem(d) and (\QPOSfive),\quad ${\Phi_A}^@ \ge f^@ \ge {\Psi_A}^@ = \Phi_A \ge q$ on $B$.   Thus, using Lemma \THPSlem(f),
$$A \subset \PQ\big({\Phi_A}^@\big) \subset \PQ\big(f^@\big) \subset \PQ(\Phi_A).\meqno\HHATthree$$
The result now follows from (\HHATtwo), (\HHATthree) and Theorem \PHIMAXthm.\qed
\defSection\HOMsec
\medbreak
\leftline{\bf \HOMsec.\quad SSD--homomorphisms and the Gossez extension}
\medskip
\noindent
The main result in this section is Theorem \AGOSSthm, in which we extend to SSD spaces a concept originally due to Gossez for maximally monotone multifunctions.
\defDefinition \DUALITYdef
\medbreak
\noindent
{\bf Definition \DUALITYdef.}\enspace Let $\big(B,\xbra\cdot\cdot\big)$ and $\big(D,\ybra\cdot\cdot\big)$ be SSD spaces.   In this case, we will \slant always\endslant\ write\quad $\qt(d) := \half\ybra{d}{d}$\quad($d \in D$).\quad We say that $\iota \colon\ B \to D$ is a \slant SSD--homomorphism\endslant\ if $\iota$ is linear and,
$$\all\ b,c\in B,\qquad \Ybra{\iota(b)}{\iota(c)} = \xbra{b}{c},\meqno\DUALone$$
from which
$$\qt \circ \iota = q.\meqno\DUALtwo$$
Let $\iota\colon B \to D$ be an SSD--homomorhism.   We define the bilinear map $\bra\cdot\cdot_\iota\colon B \times D \to \r$ by
$$\bra{b}{d}_\iota := \ybra{\iota(b)}{d}\qquad\big((b,d) \in B \times D\big).\meqno\DUALthree$$
We note then from (\DUALthree) and (\DUALone) that if $b,c \in B$ then $\Bra{b}{\iota(c)}_\iota = \Ybra{\iota(b)}{\iota(c)} = \xbra{b}{c}$, and so $\big(D,\iota,\bra\cdot\cdot_\iota\big)$ is a linked external space (see Definition \THETAdef).   If $C$ is a nonempty $\qt$--positive subset of $D$ then, as in (\QPOStwo), (\QPOSthree) and (\QPOSfour),
$$\Phi_C(d) := \supn_C\big[\ybra\cdot{d} - \qt\big] = q(d) - \inf \qt(C - d)\quad(d \in D)\quand \Phi_C = \qt\ \on\ C.\meqno\DUALfour$$
If $f \in \PC(D)$ and $d \in D$ then, as in (\FATone), 
$$f^@(d) := \supn_D\big[\ybra\cdot{d} - f\big].\meqno\DUALfive$$
The next two lemmas contain the preliminary results that we will need. We obtain Lemma \LDlem\ by transcribing Lemmas \Llem, \FFATlem(b) and \PHlem(b,c) to our present situation.
\defLemma \LDlem
\medbreak
\noindent
{\bf Lemma \LDlem.}{\rm(a)}\enspace\slant Let $f \in \PC(D)$,\quad $f \ge \qt$ on $D$\quad and $\PQT(f) \ne \emptyset$.   Then $\PQT(f)$ is a $\qt$--positive subset of $D$ and\qquad $f^@ = \qt$ on $\PQT(f)$.
\smallbreak
\noindent
{\rm(b)}\enspace If $C$ is a nonempty $\qt$--positive subset of $D$ then\quad${\Phi_C}^@ \ge  \Phi_C$\quad and\quad${\Phi_C}^{@@} = \Phi_C$ on $D$.
\endslant
\defLemma \HOMlem
\medbreak
\noindent
{\bf Lemma \HOMlem.}\enspace\slant Let $\big(B,\xbra\cdot\cdot\big)$ and $\big(D,\ybra\cdot\cdot\big)$ be SSD spaces, $\iota \colon\ B \to D$ be an SSD--\break homomorphism and $A$ be a nonempty $q$--positive subset of $B$.   Then:
\smallskip\noindent
{\rm(a)}\enspace $\iota(A)$ is a nonempty $\qt$--positive subset of $D$ and $\Phi_{\iota(A)} \in \PC(D)$.
\smallskip
\noindent
{\rm(b)}\enspace For all $d \in D$, $\Theta_A(d) = \Phi_{\iota(A)}(d) = \qt(d) - \inf\qt\big(\iota(A) - d\big)$.
\smallskip
\noindent
{\rm(c)}\enspace ${\Theta_A}^@ \ge {\Phi_A}^* \ge \Theta_A$ on $D$\quand${\Theta_A}^@ \ge {\Phi_A}^{*@} \ge \Theta_A$ on $D$.\endslant
\Proof If $b,c \in A$ then (\DUALtwo) gives $\qt\big(\iota(b) - \iota(c)\big) = \qt \circ \iota(b - c) = q(b - c) \ge 0$, and (a) follows from (\DUALfour).   If $d \in D$ then the first expression in (\LINKfour), (\DUALthree) and (\DUALtwo) give $\Theta_A(d) = \supn_A\big[\bra\cdot{d}_\iota - q\big] = \supn_A\big[\Ybra{\iota(\cdot)}{d} - \qt\circ \iota\big] = \supn_{\iota(A)}\big[\ybra{\cdot}{d} - \qt\big]$, and (b) follows from (\DUALfour). 
\smallskip
Of course, it is immediate from (b) and Lemma \LDlem(b) that\quad ${\Theta_A}^@ = {\Phi_{\iota(A)}}^@ \ge \Phi_{\iota(A)} = \Theta_A$ on $D$,\quad but for (c) we need the extra information about ${\Phi_A}^*$ and ${\Phi_A}^{*@}$.   So let $d \in D$.   Then, from (\DUALfive), (\DUALthree), (\LINKfive), (\LINKtwo) and the second expression in (\LINKfour),   
$$\eqalign{{\Theta_A}^@(d)
&= \supn_D\big[\ybra{\cdot}{d} - \Theta_A\big]
\ge \supn_B\big[\Ybra{\iota(\cdot)}{d} - \Theta_A \circ \iota\big] = \supn_B\big[\bra{\cdot}{d}_\iota - \Phi_A\big]\cr
&= {\Phi_A}^*(d) = \supn_B\big[\bra{\cdot}{d}_\iota - \Phi_A\big]
\ge \supn_A\big[\bra{\cdot}{d}_\iota - \Phi_A\big] = \Theta_A(d).}$$
This completes the proof of the first assertion in (c).   It is clear from this, (b), and Lemma \LDlem(b) that\quad ${\Theta_A}^@ \ge {\Phi_A}^{*@} \ge {\Theta_A}^{@@} = {\Phi_{\iota(A)}}^{@@} = \Phi_{\iota(A)} = \Theta_A$ on $D$,\quad which gives the second assertion in (c).\qed
\medbreak
The next concept is a generalization to SSD spaces of an idea originally introduced by Gossez in \cite\GOSSEZ\endcite\ for maximally monotone multifunctions.   The use of the word ``extension'' in Definition \GOSSEZdef\ will be justified in Theorem \AGOSSthm(a).
\defDefinition \GOSSEZdef
\medbreak
\noindent
{\bf Definition \GOSSEZdef.}\enspace Let $\big(B,\xbra\cdot\cdot\big)$ and $\big(D,\ybra\cdot\cdot\big)$ be SSD spaces, $\iota \colon\ B \to D$ be an SSD--homomorphism and $A$ be a nonempty $q$--positive subset of $B$.   It is clear from Lemma \HOMlem(b) that if $d \in D$ then\qquad $\Theta_A(d) \le \qt(d) \iff \inf\qt\big(\iota(A) - d\big) \ge 0$.  We define the \slant Gossez extension\endslant\ of $A$ in $D$ to be the set
$$\AGOSS = \big\{d \in D\colon\ \inf\qt\big(\iota(A) - d\big) \ge 0\big\} = \big\{d \in D\colon\ \Theta_A(d) \le \qt(d)\big\}.\meqno\GOSSEZone$$
\par
Theorem \AGOSSthm\ will be used in Theorem \EQthm, which will be used, in turn, in Theorem \EEEQthm, where we prove that for a maximally monotone set,  type (ED), dense type, type (D), type (WD) and type (NI) are all equivalent.   Theorem \EEEQthm\ will be used in Theorems \BRthm\ and \DFPthm.  
\defTheorem \AGOSSthm
\medbreak
\noindent
{\bf Theorem \AGOSSthm.}\enspace\slant Let $\big(B,\xbra\cdot\cdot\big)$ and $\big(D,\ybra\cdot\cdot\big)$ be SSD spaces, $\iota \colon\ B \to D$ be an SSD--homomorphism and $A$ be a nonempty $q$--positive subset of $B$.   Then:
\smallbreak\noindent
{\rm(a)}\enspace $\iota(A) \subset \AGOSS$.
\smallbreak\noindent
{\rm(b)}\enspace If, for all $d \in \AGOSS$, $\inf\qt\big(\iota(A) - d\big) \le 0$ then
$$\Theta_A \ge \qt\ \on\ D.\meqno\AGOSSone$$
{\rm(c)}\enspace If {\rm(\AGOSSone)} is satisfied then $\AGOSS =\PQT({\Phi_A}^{*@}) = \PQT({\Phi_A}^*) = \PQT({\Theta_A}^@) = \PQT(\Theta_A)$.
\endslant
\Proof If $a \in A$ then, from  (\LINKfive), (\QPOSfour) and  (\DUALtwo), $\Theta_A\big(\iota(a)\big) = \Phi_A(a) = q(a) = \qt\circ \iota(a) = \qt\big(\iota(a)\big)$, and so the second expression in (\GOSSEZone) implies that $\iota(a) \in \AGOSS$.   This gives (a).   (b) is immediate from Lemma \HOMlem(b).   Suppose, finally, that (\AGOSSone) is satisfied.   It is obvious from (a) and the second expression in (\GOSSEZone) that $\emptyset \ne \AGOSS = \PQT\big(\Theta_A\big)$ and, from Lemma \HOMlem(c), that
$$\PQT({\Theta_A}^@) \subset \PQT({\Phi_A}^*) \subset \PQT(\Theta_A) \quand \PQT({\Theta_A}^@) \subset \PQT({\Phi_A}^{*@}) \subset \PQT(\Theta_A).$$
(c) follows since Lemma \LDlem(a) with $f := \Theta_A $ gives $\PQT(\Theta_A) \subset \PQT({\Theta_A}^@)$.\qed
\defSection\BSSDsec
\medbreak
\leftline{\bf \BSSDsec.\quad VZ functions on Banach SSD spaces}
\medskip
\noindent
We note that we do not use anything from Section \HOMsec\ in this section --- we will combine this section with Section \HOMsec\ in Section \DBSSDsec.   The other comment is that the dual of a Banach space is not mentioned explicitly until Remark \BSTARrem, though there is an implicit use of the dual in Theorem \VZEXthm(b) (in the observation that a proper convex lower semicontinuous function dominates a continuous affine function).   If $X$ is a nonzero real Banach space, we write $\PCLSC(X)$ for the set 
$$\{f \in \PC(X)\colon\ f\ \hbox{is lower semicontinuous on}\ X\}.$$\par
%
\defDefinition \NORMdef
\noindent
{\bf Definition \NORMdef.}\enspace We say that $\big(B,\xbra\cdot\cdot,\|\cdot\|\big)$ is a \slant Banach SSD space\endslant\ if $\big(B,\xbra\cdot\cdot\big)$ is an SSD space and $\|\cdot\|$ is a norm on $B$ with respect to which $B$ is a Banach space and
$$\all\ b,c\in B,\quad \big|\xbra{b}{c}\big| \le \|b\|\|c\|.\meqno\IOTAtwo$$
If we take $c = b$, we derive that
$$\half\|\cdot\|^2 + q \ge 0\ \on\ B.\meqno\IOTAone$$
Then, for all $d,e \in B$,
$$|q(d) - q(e)| = \half\big|\xbra{d}{d} - \xbra{e}{e}\big| = \half\big|\xbra{d- e}{d + e}\big| \le \half\|d - e\|\|d + e\|.\meqno\IOTAthree$$
We define the continuous even functions $g$ and $p$ on $B$ by $g := \half\|\cdot\|^2$ and $p := g + q$, so that\quad $p \ge 0$ on $B$.\quad Since $p(0) = 0$, in fact
$$\infn_Bp = 0.\meqno\IOTAfour$$
Also, for all $d,e \in B$, $|g(d) - g(e)| = \half\big|\|d\| - \|e\|\big|\big(\|d\| + \|e\|\big) \le \half\|d - e\|\big(\|d\| + \|e\|\big)$.   Combining this with (\IOTAthree), for all $d,e \in B$,
$$|p(d) - p(e)| \le \|d - e\|\big(\|d\| + \|e\|\big).\meqno\IOTAfive$$
\par
\defRemark \NORMrem
\medbreak
\noindent
{\bf Remark \NORMrem.}\enspace  It is clear from (\IOTAtwo) that Example \Hex\ is a Banach SSD space if, and only if $\|T\| \le 1$, which is the case with (a), (b) and (c).   Looking ahead to Remark \BSTARrem, it then follows from (\IOTAnew) that $\iota = T$.  
\defExample \EENex
\medbreak
\noindent
{\bf Example \EENex.}\enspace We now continue our discussion of Example \EEex.   We suppose that $E$ is a nonzero Banach space, $B = E \times E^*$ and, for all $(x,x^*) \in B$, $\|(x,x^*)\|_2 := \sqrt{\|x\|^2 + \|x^*\|^2}$.   It is clear from the CauchyÐ-Schwarz inequality that (\IOTAtwo) is satisfied, and so\break $\big(E \times E^*,\xbra\cdot\cdot,\|\cdot\|_2\big)$ is a Banach SSD space.
\medbreak
We now introduce the concept of inf--convolution.   This certainly goes back as far as \cite\RTRCA\endcite, but we emphasize that we will be using it here for nonconvex functions.     
\defDefinition \VZdef
\medbreak
\noindent
{\bf Definition \VZdef.}\enspace
Let $X$ be a vector space and $h,k\colon X \to \rbar$.   The \slant inf--convolution of $h$ and $k$\endslant\ is defined by\quad $(h \episum k)(x) := \infn_{y \in X} \big[h(y) + k(x - y)\big]$\quad($x \in X$).\quad   It is clear that
$$\infn_Xk = 0 \qlr \infn_X\big[h \episum k\big] = \infn_Xh.\meqno\VZone$$
Now let $\big(B,\xbra\cdot\cdot,\|\cdot\|\big)$ be a Banach SSD space and $f \in \PC(B)$.   We say that $f$ is a \slant VZ function\endslant\ if
$$(f - q) \episum p = 0\ \on\ B.\meqno\VZtwo$$
It follows from (\IOTAfour) and (\VZone) that
$$\hbox{if}\ f\ \hbox{is a VZ function then}\ \infn_B[f - q] = 0.\meqno\VZthree$$
``VZ'' stands for ``Voisei--Z\u alinescu'', since (\VZtwo) is an extension to Banach SSD spaces of a condition introduced in \cite\VZ, Proposition 3\endcite.   The following simple inequality will be useful:  suppose that $f \in \PC(B)$ and\quad $f \ge q$ on $B$;\quad then, for all $c \in B$,
$$\left.\eqalign{\big((f - q) \episum p\big)(c) &= \infn_{b \in B}\big[(f - q)(b) + p(c - b)\big]\cr
&\le \infn_{a \in \PQ(f)}\big[(f - q)(a) + p(c - a)\big]\cr
&= \inf p\big(c - \PQ(f)\big) = \inf p\big(\PQ(f) - c\big).}\right\}\meqno\VZfour$$
%
\defDefinition \Pdef
\medbreak
\noindent
{\bf Definition \Pdef.}\enspace Let $A$ be a subset of a Banach SSD space $\big(B,\xbra\cdot\cdot,\|\cdot\|\big)$.   We say that $A$ is \slant $p$--dense in\endslant\ $B$ if, for all $c \in B$, $\inf p(A - c) = 0$.     
\medbreak
We now come to our main results on VZ functions on Banach SSD spaces.   We shall see in Remark \FPHIrem\ that the constant $\sqrt 2$ in (\VZEXtwo) is sharp, and also that (\VZEXtwo) leads to a strict strengthening of (\VZEXone).   If we take Theorem \VZEXthm(a) into account then a double induction is used to prove Theorem \VZEXthm(c).   We do not know if this is actually necessary.    Theorem \VZEXthm(d) is an extension to Banach SSD spaces of \cite\VZ, Theorem 8\endcite.    Theorem  \VZEXthm\ is related to some results proved by Zagrodny in \cite\ZAGRODNY\endcite.   These are discussed more fully in the comments preceding Problem \ZAGprob. 
\defTheorem  \VZEXthm
\medbreak
\noindent
{\bf Theorem  \VZEXthm.}\enspace\slant Let $\big(B,\xbra\cdot\cdot,\|\cdot\|\big)$ be a Banach SSD space and $f \in \PCLSC(B)$ be a VZ function.   Then:
\smallbreak\noindent
{\rm(a)}\enspace $\PQ(f)$ is a $q$--positive subset of $B$ and
$$d \in \dom\,f \qlr \dist(d,\PQ(f)) \le 5\sqrt{(f - q)(d)}.\meqno\VZEXone$$
{\rm(b)}\enspace $\PQ(f)$ is $p$--dense in $B$.
\smallbreak\noindent
{\rm(c)}\enspace For all $c \in B$,  $\inf q\big(\PQ(f) - c\big) \le 0$ and
$$c \in B \qlr \dist\big(c,\PQ(f)\big) \le \rttwo\sqrt{{-}\inf q\big(\PQ(f) - c\big)}.\meqno\VZEXtwo$$
{\rm(d)}\enspace  $\PQ(f)$ is a maximally $q$--positive subset of $B$.
\endslant
\Proof (a)\enspace (\VZthree) implies that\quad $f \ge q$ on $B$,\quad and so $\PQ(f)$ is defined.   Let $d \in \dom\,f$.   We first prove that there exists a Cauchy sequence $\{b_n\}_{n \ge 1}$ such that, for all $n \ge 1$,
$$(f - q)(b_n) \le (f - q)(d)/4^n \quand \|d - b_n\| \le 5\sqrt{(f - q)(d)}.\meqno\EXtwo$$
Since we can take $b_n = d$ if $(f - q)(d) = 0$, we can and will suppose that
$$\alpha := \sqrt{(f - q)(d)} > 0.\meqno\EXthree$$
Write $b_0 := d$.   Then we can choose inductively $b_1,b_2,\dots \in B$ \big(using the fact that\break $\big((f - q) \episum p\big)(b_{n - 1}) = 0$\big) such that, for all $n \ge 1$,\quad $(f - q)(b_n) + p(b_{n - 1} - b_n) \le (\alpha/2^n)^2$.\quad It follows from (\VZthree), (\EXthree) and (\IOTAfour) that, 
$$\all\ n \ge 1,\qquad p(b_{n - 1} - b_n) \le (\alpha/2^n)^2,\meqno\EXfour$$
and
$$\all\ n \ge 0,\qquad (f - q)(b_n) \le (\alpha/2^n)^2.\meqno\EXfive$$
Substituting this into Lemma \ROOTlem, for all $n \ge 1$,
$$-q(b_{n - 1} - b_n) \le \Big[\sqrt{(f - q)(b_{n - 1})} + \sqrt{(f - q)(b_n)}\Big]^2 \le [\alpha/2^{n - 1} + \alpha/2^n]^2 = 9(\alpha/2^n)^2.$$
Consequently, since $g(b_{n - 1} - b_n) = p(b_{n - 1} - b_n) - q(b_{n - 1} - b_n)$, (\EXfour) gives,
$$\all\ n \ge 1,\qquad g(b_{n - 1} - b_n) \le (\alpha/2^n)^2 + 9(\alpha/2^n)^2 = 10(\alpha/2^n)^2,$$
and so, for all $n \ge 1$, $\|b_{n - 1} - b_n\| \le 5\alpha/2^n$.   Adding up this inequality for $n = 1, \dots, m$ and using (\EXthree), we derive that, for all $m \ge 1$, $\|d - b_m\| \le 5\alpha = 5\sqrt{(f - q)(d)}$, and (\EXtwo) now follows from (\EXfive).   Now set $a = \lim_{n}b_n$, so that $\|d - a\| \le 5\sqrt{(f - q)(d)}$.   (\EXtwo) and the lower semicontinuity of $f - q$ now imply that $(f - q)(a) \le 0$, that is to say, $a \in \PQ(f)$.   Since $\dom\,f \ne \emptyset$, it follows that $\PQ(f) \ne \emptyset$ and so, from Lemma \Llem, $\PQ(f)$ is a $q$--positive subset of $B$, and obviously (\VZEXone) is satisfied.   This completes the proof of (a).
\smallskip
(b)\enspace Let $c \in B$.   Since $\big((f - q) \episum p\big)(c) = 0$, we can choose inductively $d_1,d_2,\dots \in B$ such that, for all $n \ge 1$,
$$f(d_n) + g(c - d_n) + q(c) - \xbra{c}{d_n} = (f - q)(d_n) + p(c - d_n) < 1/n^2.$$
Consequently, using (\IOTAfour), (\VZthree) and (\IOTAtwo), for all $n \ge 1$,
$$(f - q)(d_n) <  1/n^2,\ p(c - d_n) < 1/n^2\meqno\EXseven$$
and
$$f(d_n) + g(c - d_n) + q(c) - \|c\|\|d_n\| < 1/n^2.\meqno\EXeight$$
Since $f \in \PCLSC(B)$, $f$ dominates a continuous affine function, and so (\EXeight) and the usual coercivity argument imply that $K := \sup_{n \ge 1}\|d_n\| < \infty$.   From (a) and (\EXseven), for all $n \ge 1$, there exists $a_n \in \PQ(f)$ such that $\|a_n - d_n\| \le 5/n$.   Now, from (\IOTAfive), for all $n \ge 1$,
$$\eqalign{|p(c - a_n) - p(c - d_n)|
&\le \|a_n - d_n\|(2\|c\| + \|a_n\| + \|d_n\|)\cr
&\le \big(2\|c\| + (K + 5) + K\big)5/n.}$$
Thus $\lim_{n \to \infty}\big[p(c - a_n) - p(c - d_n)\big] = 0$, and (b) follows by combining this with (\EXseven).
\smallbreak
(c)\enspace Let $c \in B$.   Then, from (b),
$$\inf g\big(\PQ(f) - c\big) + \inf q\big(\PQ(f) - c\big) \le \inf (g + q)\big(\PQ(f) - c\big) = \inf p\big(\PQ(f) - c\big) = 0.$$
Thus
$\half\dist\big(c,\PQ(f)\big)^2 = \inf g\big(\PQ(f) - c\big) \le - \inf q\big(\PQ(f) - c\big)$, from which (\VZEXtwo) is an immediate consequence.
\smallskip
(d)\enspace We suppose that $c \in B$ and $\inf q\big(\PQ(f) - c\big) \ge 0$, and we must prove that $c \in \PQ(f)$.   From (c), in fact\quad $\inf q\big(\PQ(f) - c\big) = 0$ \quand $\dist\big(c,\PQ(f)\big) = 0$.\quad   The lower semicontinuity of $f$ implies that $\PQ(f)$ is closed, and so $c \in \PQ(f)$.   This completes the proof of (d).\qed
\medbreak
The interest of Theorem \VZEQthm\ below is that it tells us that we can determine whether $f$ is a VZ function by simply inspecting $\PQ(f)$.
\defTheorem \VZEQthm
\medbreak
\noindent
{\bf Theorem \VZEQthm.}\enspace\slant Let $\big(B,\xbra\cdot\cdot,\|\cdot\|\big)$ be a Banach SSD space and $f \in \PCLSC(B)$.   Then\quad $f$ is a VZ function $\quad\iff\quad$ $f \ge q$ on $B$\quad and\enspace $\PQ(f)$ is $p$--dense in $B$.\endslant
\Proof ($\lr$) is immediate from (\VZthree) and Theorem \VZEXthm(b). Suppose, conversely, that $f \ge q$ on $B$\quad and\enspace $\PQ(f)$ is $p$--dense in $B$.   Then from (\VZfour), for all $c \in B$, 
$\big((f - q) \episum p\big)(c) \le \inf p\big(\PQ(f) - c\big) = 0$, from which\quad $(f - q) \episum p \le 0$ on $B$.\quad   On the other hand, since\quad $f - q \ge 0$ on $B$\quad and, from (\IOTAfour),\quad $p \ge 0$ on $B$,\quad we have\quad $(f - q) \episum p \ge 0$ on $B$.\quad   Thus $f$ is a VZ function, as required.\qed
\medbreak
We point out that the function $h$ in Theorem \EXthm(a) is not required to be lower semicontinuous, so we cannot simply apply Theorem \VZEQthm\ with $f$ replaced by $h$. 
\defTheorem \EXthm
\medbreak
\noindent
{\bf Theorem \EXthm.}\enspace\slant Let $\big(B,\xbra\cdot\cdot,\|\cdot\|\big)$. be a Banach SSD space and $f \in \PCLSC(B)$ be a VZ function.   Then: 
\smallskip\noindent
{\rm(a)}\enspace Let $h \in \PC(B)$, $h \ge q$ on $B$, and $\PQ(h) \supset \PQ(f)$.   Then $\PQ(h) = \PQ(f)$ and $h$ is a VZ function.
\smallskip\noindent
{\rm(b)}\enspace $f^@ \in \PCLSC(B)$, $f^@$ is a VZ function and $\PQ\big(f^@\big) = \PQ(f)$.\endslant
\Proof(a)\enspace It is clear from Theorem \VZEXthm(d) that $\PQ(h) = \PQ(f)$.   From (\VZfour) and Theorem \VZEXthm(b), for all $c \in B$, $\big((h - q) \episum p\big)(c) \le \inf p\big(\PQ(h) - c\big) = \inf p\big(\PQ(f) - c\big) = 0$,
from which\quad $(h - q) \episum p \le 0$ on $B$.\quad   On the other hand, since\quad $h - q \ge 0$ on $B$\quad and\quad $p \ge 0$ on $B$,\quad we have\quad $(h - q) \episum p \ge 0$ on $B$.\quad   Thus $h$ is a VZ function.
\smallskip
(b)\enspace Let $c \in B$. Then, since\quad $q \le p$ on $B$,\quad Definition \FATdef\ gives
$$q(c) - f^@(c) = \infn_{b \in B}\big[f(b) - \xbra{b}{c} + q(c)\big] = \big((f - q) \episum q\big)(c) \le \big((f - q) \episum p\big)(c)= 0,$$
and so\quad $f^@ \ge q$ on $B$.\quad   It now follows from Lemma \FFATlem(b) that $\PQ\big(f^@\big) \supset \PQ(f)$, and so (a) implies that $f^@$ is a VZ function and $\PQ\big(f^@\big) = \PQ(f)$.   Since $\PQ(f) \ne \emptyset$, it is evident that $f^@ \in \PCLSC(B)$.\qed
\defRemark \BSTARrem
\medbreak
\noindent
{\bf Remark \BSTARrem.}\enspace Up to this point, we have not mentioned the Banach space dual, $B^*$, of $B$.   It is easy to see from (\IOTAtwo) and standard algebraic arguments that there exists a linear map $\iota\colon\ B \to B^*$ such that $\|\iota\| \le 1$ and
$$\all\ b,c\in B,\quad \Bra{b}{\iota(c)} = \xbra{b}{c}.\meqno\IOTAnew$$
It follows that $\big(B^*,\iota,\bra\cdot\cdot\big)$ is a linked external space (see Definition \THETAdef) and, if $A$ is a nonempty $q$--positive subset of B, we define $\Theta_A \in \PC(B^*)$ and $\Psi_A \in \PC(B)$ by (\LINKfour) and (\LINKsix), with $D$ replaced by $B^*$.  
\medbreak
The proof of Theorem \EXthm\ relies heavily on the lower semicontinuity of $f$.   We will show in Corollary \Hcor\ below that part of Theorem \EXthm(b) can be recovered even if $f$ is not assumed to be lower semicontinuous.
\defCorollary \Hcor
\medbreak
\noindent
{\bf Corollary \Hcor.}\enspace\slant Let $\big(B,\xbra\cdot\cdot,\|\cdot\|\big)$ be a Banach SSD space and $f \in \PC(B)$ be a VZ function.   Then $f^@ \in \PCLSC(B)$, $f^@$ is a VZ function and $\PQ\big(f^@\big)$ is  a maximally $q$--positive subset of $B$.\endslant
\Proof Let $\fbar$ be the lower semicontinuous envelope of $f$.   Since $q$ is continuous and\quad $f \ge q$ on $B$,\quad it follows that\quad $f \ge \fbar \ge q$ on $B$.\quad   Thus, from  (\IOTAfour),
$$0 = (f - q) \episum p \ge (\fbar - q) \episum p \ge 0 \episum p = 0\ \on\ B,$$      
and so $\fbar$ is a VZ function.   Since $\fbar \in \PCLSC(B)$, Theorem \EXthm(b) implies that $\fbar^@$ is a VZ function also.   It is well known that $\fbar^* = f^*$ on $B^*$ thus, composing with $\iota$ and using (\LINKthree),\quad $\fbar^@ = f^@$ on $B$.\quad The result now follows from Theorem \VZEXthm(d), with $f$ replaced by $f^@$.\qed
\defTheorem \VZONLYthm
\medbreak
\noindent
{\bf Theorem \VZONLYthm.}\enspace\slant Let $\big(B,\xbra\cdot\cdot,\|\cdot\|\big)$ be a Banach SSD space, $f \in \PCLSC(B)$ be a VZ function    and $A := \PQ(f)$.   Then
$$\Psi_A \ge f \ge \Phi_A \ge q\ \on\ B \quand {\Phi_A}^* \ge f^* \ge \Theta_A\ \on\ B^*,\meqno\VZONLYone$$
$$\PQ\big(\Psi_A\big) = \PQ\big({\Phi_A}^@\big) = \PQ\big(\Phi_A\big) = A,\meqno\VZONLYtwo$$
and
$$\Phi_A, {\Phi_A}^@\ \and\ \Psi_A\ \hbox{are all VZ functions}.\meqno\VZONLYthree$$
Now let $h \in \PC(B)$ \quand $\Psi_A \ge h \ge \Phi_A$ on $B$.\quad   Then $h$ and $h^@$ are VZ functions.\endslant
\Proof We first note from (\VZthree) and Theorem \VZEXthm(d) that\quad $f \ge q$ on $B$\quad and $A$ is a maximally $q$--positive subset of $B$. From standard normed space theory, $f$ is $w(B,B^*)$--lower semicontinuous and so (\VZONLYone) and (\VZONLYtwo) follow from Corollary \MAXcor\ and (\QPOSfive).   The first assertion in (\VZONLYone), (\VZONLYtwo) and Theorem \VZEQthm\ imply that $\Psi_A$ and $\Phi_A$ are VZ functions, and (\VZONLYthree) now follows from Theorem \EXthm(b), with $f$ replaced by $\Phi_A$.   The assertions about $h$ and $h^@$ follow from Theorem \EXthm(a) and Corollary \Hcor.\qed
\defRemark \FPHIrem
\medbreak
\noindent
{\bf Remark \FPHIrem.}\enspace Let $\big(B,\xbra\cdot\cdot,\|\cdot\|\big)$ be a Banach SSD space and $f \in \PCLSC(B)$ be a VZ function.   We know from Theorem \VZONLYthm\ that $\PQ\big(\Phi_{\PQ(f)}\big) = \PQ(f)$, $\Phi_{\PQ(f)}$ is a VZ function and\quad $\Phi_{\PQ(f)} \le f$ on $B$.\quad Combining this with (\QPOSthree), for all $c \in B$,
$$-\inf q\big(\PQ(f) - c\big) = \big(\Phi_{\PQ(f)} - q\big)(c) \le (f - q)(c).\meqno\FPHIone$$
Thus Theorem \VZEXthm(c) implies that, for all $c \in B$,
$$\dist(c,\PQ(f)) = \dist\big(c,\PQ\big(\Phi_{\PQ(f)}\big)\big) \le \rttwo\sqrt{\big(\Phi_{\PQ(f)} - q\big)(c)} \le \rttwo\sqrt{(f - q)(c)}.\meqno\FPHItwo$$  
This shows that Theorem \VZEXthm(c) is stronger than Theorem  \VZEXthm(a).   Now consider the Banach SSD space $\big(\r \times \r,\xbra\cdot\cdot,\|\cdot\|_2\big)$, where the notation is as in Examples \EENex.   Define $f \in \PCLSC(B)$ by $f(x_1,x_2) := \half(x_1^2 + x_2^2)$.   Then $(f - q)(x_1,x_2) = \half(x_1^2 + x_2^2) - x_1x_2 = \half(x_1 - x_2)^2$
and
$p(x_1,x_2) = \half(x_1^2 + x_2^2) + x_1x_2 = \half(x_1 + x_2)^2$.
Let $c := (z_1,z_2) \in B$ and $b := \big(\half(z_1 + z_2),\half(z_1 + z_2)\big) \in B$.   Then $(f - q)(b) = 0$ and $p(c - b) = 0$.   Consequently, $f$ is a VZ function.  Now $\PQ(f)$ is the diagonal of $\r^2$ and so, by direct computation, for all $c = (x_1,x_2) \in \r^2$,
$-\inf q\big(\PQ(f) - c\big) = \fourth(x_1 - x_2)^2$.   Since $\fourth(x_1 - x_2)^2 < \half(x_1 - x_2)^2$ when $x_1 \ne x_2$, the inequality in (\FPHIone) is generally strict.\smallbreak
Now let $h := \Phi_{\PQ(f)}$.   (\QPOSthree) gives us that, for all $(x_1,x_2) \in B$,
$$\sqrt{(h - q)(x_1,x_2)} = \sqrt{\fourth(x_1 - x_2)^2} = \half|x_1 - x_2|.$$
On the other hand, by direct computation, $\dist\big((x_1,x_2),\PQ(h)\big) = \rthalf|x_1 - x_2|$.   Thus the constant $\rttwo$ in the inequalities in (\FPHItwo) is sharp.   The genesis of this argument and example can be found in the results of Mart\'\i nez-Legaz and Th\'era in \cite\MT\endcite.
\defRemark \TWOrem
\medbreak
\noindent
{\bf Remark \TWOrem.}\enspace The following result follows by applying the inequality between the first and last terms in (\FPHItwo) to Example \EENex:  \slant Let $E$ be a nonzero Banach space, and $f$ be a lower semicontinuous VZ function on $E \times E^*$.
Then, for all $(x,x^*) \in E \times E^*$,
$$\infn_{(y,y^*) \in \PQ(f)}\sqrt{\|y - x\|^2 + \|y^* - x^*\|^2}
\le \rttwo\sqrt{f(x,x^*) - \bra{x}{x^*}}.$$\endslant
This strengthens the result proved in \cite\VZ, Theorem 4\endcite, namely that
$$\infn_{(y,y^*) \in \PQ(f)}\sqrt{\|y - x\|^2 + \|y^* - x^*\|^2} \le 2\sqrt{f(x,x^*) - \bra{x}{x^*}}.$$
As we observed in Remark \FPHIrem, the constant $\sqrt 2$ is sharp.
\defRemark \VZONLYrem
\medbreak
\noindent
{\bf Remark \VZONLYrem.}\enspace We note that the inequalities for $B$ in (\VZONLYone) have four functions, while the inequality for $B^*$ has only three.   The reason for this is that we do not have a function on $B^*$ that plays the role that the function $q$ plays on $B$.   The function $\qt$, which plays such a role, will be introduced in this context in Definition \BSSDdef.
\defSection\DBSSDsec
\medbreak
\leftline{\bf\DBSSDsec.\quad Banach SSD duals}
\medskip
\noindent
If $X$ is a nonzero real Banach space, we write $X\dbs$ for the bidual of $X$ \big(with the pairing $\bra\cdot\cdot\colon X^* \times X\dbs \to \r$\big).   If $x \in X$, we write $\wh x$ for the canonical image of $x$ in $X\dbs$, that is to say
$$x \in X\ \and\ x^* \in X^* \qlr\bra{x^*}{\wh x} = \bra{x}{x^*}.$$\par
\defDefinition \BSSDdef
\medbreak
\noindent
{\bf Definition \BSSDdef.}\enspace Let $\big(B,\xbra\cdot\cdot,\|\cdot\|\big)$ be a Banach SSD space, $(B^*,\|\cdot\|)$ be the Banach space dual of $B$ and the bounded linear map $\iota\colon\ B \to B^*$ be defined as in (\IOTAnew).   Let $\big(B^*,\ybra\cdot\cdot,\|\cdot\|\big)$ also be a Banach SSD space.   We say that $\big(B^*,\ybra\cdot\cdot,\|\cdot\|\big)$ is a \slant Banach SSD dual\endslant\ of $\big(B,\xbra\cdot\cdot,\|\cdot\|\big)$ if $\bra\cdot\cdot_\iota = \bra\cdot\cdot$ on $B \times B^*$ \big(see (\DUALthree)\big), that is to say 
$$\all\ b \in B\ \and\ c^* \in B^*,\qquad \Ybra{\iota(b)}{c^*} = \bra{b}{c^*}.\meqno\BSSDone$$
We have not required explicitly that $\iota$ be an SSD-homomorphism from $\big(B,\xbra\cdot\cdot\big)$ into $\big(B^*,\ybra\cdot\cdot\big)$ \big(see  (\DUALone)\big):  this is automatically satisfied since (\BSSDone) and (\IOTAnew) imply that, for all $b,c \in B$, $\Ybra{\iota(b)}{\iota(c)} = \Bra{b}{\iota(c)} = \xbra{b}{c}$.
\smallskip
Thus if $\big(B,\xbra\cdot\cdot,\|\cdot\|\big)$ is a Banach SSD space with Banach SSD dual $\big(B^*,\ybra\cdot\cdot,\|\cdot\|\big)$, we can use all the results of Sections \EXTsec\ and \HOMsec\ (with ``$D$'' replaced by ``$B^*$'') and Section \BSSDsec.
\smallskip
By analogy with (\IOTAnew), we define the bounded linear map $\it\colon B^* \to B\dbs$ so that
$$\all\ c^*, b^* \in B^*,\qquad \Bra{c^*}{\it(b^*)} = \ybra{c^*}{b^*}.\meqno\BSSDthree$$
and the function $\pt\colon\ B^* \to \r$ by $\pt := \half\|\cdot\|^2 + \qt$.   Thus we have
$$\pt \ge 0\ \on\ B^*.\meqno\BSSDfour$$
\par
We now show that Definition \BSSDdef\ also leads to an automatic factorization of the canonical map from $B$ into $B\dbs$.    Lemma \FACTORlem\ will be used in Lemma \TWOTOPlem. 
\defLemma \FACTORlem
\medbreak
\noindent
{\bf Lemma \FACTORlem.}\enspace\slant Let $\big(B,\xbra\cdot\cdot,\|\cdot\|\big)$ be a Banach SSD space with Banach SSD dual\break $\big(B^*,\ybra\cdot\cdot,\|\cdot\|\big)$.   Then, for all $b \in B$, $\wh b = \it\circ\iota(b)$.\endslant
\Proof Let $b \in B$ and $c^* \in B^*$.   Then, from the definition of $\wh b$, (\BSSDone) and (\BSSDthree),
$$\Bra{c^*}{\wh b} = \bra{b}{c^*} = \Ybra{\iota(b)}{c^*} = \Ybra{c^*}{\iota(b)} = \Bra{c^*}{\it\circ\iota(b)}.$$
This gives the required result.\qed
\defRemark \Hrem
\medbreak
\noindent
{\bf Remark \Hrem.}\enspace In this remark, we suppose that the notation is as in Example \Hex, and that $\big(B,\xbra\cdot\cdot,\|\cdot\|\big)$ is a Banach SSD space with Banach SSD dual $\big(B,\ybra\cdot\cdot,\|\cdot\|\big)$.  We shall show  that $\xbra\cdot\cdot = \ybra\cdot\cdot$ on $B \times B$.
\smallbreak
We know already from Remarks \NORMrem\ and \BSTARrem\ that $\iota = T$ and  $\|\iota\| \le 1$.  We write $I_B$ for the identity map on $B$.
\smallbreak
It is clear from Lemma \FACTORlem\ that $\it\circ \iota = I_B$.   Now from (\BSSDthree), for all $b,c \in B$, $\Bra{\it(b)}{c} = \Bra{c}{\it(b)} = \ybra{c}{b}= \ybra{b}{c} = \Bra{b}{\it(c)}$, so $\it$ is self--adjoint.   For all $b \in B$, we have\break $\|\it(b) - \iota(b)\|^2 = \Bra{\it(b)}{\it(b)} - 2\Bra{\it(b)}{\iota(b)} + \bra{\iota(b)}{\iota(b)} = \big\|\it(b)\big\|^2 - 2\Bra{b}{\it \circ \iota(b)} +  \|\iota(b)\|^2 = \big\|\it(b)\big\|^2 - 2\Bra{b}{b} + \|\iota(b)\|^2 = \big\|\it(b)\big\|^2 - 2\|b\|^2 + \|\iota(b)\|^2$.   Since $\|\iota\| \le 1$ and $\|\it\| \le 1$, $\|\it(b) - \iota(b)\|^2 \le 0$, from which $\it(b) = \iota(b)$. Thus $\it = \iota$, from which $\ybra\cdot\cdot = \xbra\cdot\cdot$ as required.   In other words, $\big(B,\xbra\cdot\cdot,\|\cdot\|\big)$ is its own Banach SSD dual. 
\medbreak
The following concept will be critical in Theorem \MASVZthm. 
\defDefinition \PSEUDOdef
\medbreak
\noindent
{\bf Definition \PSEUDOdef.}\enspace Let $\big(B,\xbra\cdot\cdot,\|\cdot\|\big)$ be a Banach SSD space with Banach SSD dual\break $\big(B^*,\ybra\cdot\cdot,\|\cdot\|\big)$.
In line with Definition \Pdef, we say that $\iota(B)$ is \slant $\pt$--dense in\endslant\ $B^*$ if
$$\all\ b^* \in B^*,\quad \inf\pt\big(\iota(B) - b^*\big) = 0.\meqno\PSEUDOone$$\par
\defExample \EEDUALex
\medbreak
\noindent
{\bf Example \EEDUALex.}\enspace We now show that the Banach SSD space $\big(B,\xbra\cdot\cdot,\|\cdot\|_2\big)$ of Example \EENex\ has a Banach SSD dual, and
$$\iota(B)\ \hbox{is}\ \pt\hbox{--dense in}\ B^*.\meqno\EEfive$$
(See Definition \PSEUDOdef.)    We recall that $E$ is a nonzero Banach space and $B = E \times E^*$.   We represent $B^*$ by $E^* \times E\dbs$, under the pairing $\Bra{(x,x^*)}{(y^*,y\dbs)} := \bra{x}{y^*} + \bra{x^*}{y\dbs}$.   Then, for all $(x,x^*) \in B$, $q(x,x^*) = \bra{x}{x^*}$ and, from (\IOTAnew),  $\iota(x,x^*) = (x^*,\wh{x})$.   The dual norm on $E^* \times E\dbs$ is given by  $\|(y^*,y\dbs)\|_2 := \sqrt{\|y^*\|^2 + \|y\dbs\|^2}$.
\smallskip
Replacing $E$ by $E^*$ in Example \EENex, we define the symmetric bilinear form\break $\ybra{\cdot}{\cdot}\colon\ B^* \times B^* \to \r$ by $\Ybra{(x^*,x\dbs)}{(y^*,y\dbs)} := \bra{y^*}{x\dbs} + \bra{x^*}{y\dbs}$.   Then $\big(B^*,\ybra\cdot\cdot,\|\cdot\|_2\big)$ is a Banach SSD space.   We represent $B\dbs = (B^*)^*$ by $E\dbs \times E^{***}$ under the pairing $\Bra{(y^*,y\dbs)}{(w\dbs,w^{***})} := \bra{y^*}{w\dbs} + \bra{y\dbs}{w^{***}}$.   Then, for all $(y^*,y\dbs) \in B^*$, $\qt(y^*,y\dbs) = \bra{y^*}{y\dbs}$ and $\it(y^*,y\dbs) = (y\dbs,\wh{y^*})$.
\smallskip
We next show that $\big(B^*,\ybra\cdot\cdot,\|\cdot\|_2\big)$ is a Banach SSD dual of $\big(B,\xbra\cdot\cdot,\|\cdot\|_2\big)$, that is to say (\BSSDone) is satisfied.   To this end, let $(x,x^*) \in B$ and $(y^*,y\dbs) \in B^*$.   Then $\Ybra{\iota(x,x^*)}{(y^*,y\dbs)} = \Ybra{(x^*,\wh{x})}{(y^*,y\dbs)} = \bra{y^*}{\wh x} + \bra{x^*}{y\dbs} = \bra{x}{y^*} + \bra{x^*}{y\dbs} = \Bra{(x,x^*)}{(y^*,y\dbs)}$, which gives (\BSSDone), as required.
\smallskip
We now establish (\EEfive).   To see this, let $(y^*,y\dbs) \in B^*$ and $\eps > 0$.   From the definition of $\|y\dbs\|$, there exists $z^* \in E^*$ such that $\|z^*\| \le \|y\dbs\|$ and $\bra{z^*}{y\dbs} \ge \|y\dbs\|^2 - \eps$.   But then
$$\iota(0,y^* + z^*) - (y^*,y\dbs) = (y^* + z^*,0) - (y^*,y\dbs) = (z^*,-y\dbs) \in B^*.$$
Since $$\pt(z^*,-y\dbs) =  \half\big(\|z^*\|^2 + \|y\dbs\|\big) - \bra{z^*}{y\dbs} \le \|y\dbs\|^2 - \bra{z^*}{y\dbs} \le \eps,$$ we have established (\EEfive), as required.
\medbreak
The following observation will be useful in our discussion of monotone sets in Section \MONsec:  if $(a,a^*) \in B$ and $(y^*,y\dbs) \in B^*$ then\quad $\qt\big(\iota(a,a^*) - (y^*,y\dbs)\big) = \qt\big((a^*,\wh{a}) - (y^*,y\dbs)\big) = \qt(a^* - y^*,\wh{a} - y\dbs) = \bra{a^* - y^*}{\wh{a} - y\dbs}$.\quad   As a consequence, if $\emptyset \ne A \subset B$ then
$$\inf\qt\big(\iota(A) - (y^*,y\dbs)\big) =  \infn_{(a,a^*) \in A}\bra{a^* - y^*}{\wh a - y\dbs}.\meqno\EEeight$$\par
\defRemark \OTHERrem
\medbreak
\noindent
{\bf Remark \OTHERrem.}\enspace In the situation of Example \EEDUALex, there are norms $\|\cdot\|$ on $B$ other than $\|\cdot\|_2$ under which $\big(B,\xbra\cdot\cdot,\|\cdot\|\big)$ has a Banach SSD dual, and $\iota(B)$ is $\pt$--dense in $B^*$.   We refer the reader to \cite\NRSSD, Example 2.4, p. 6\endcite\ and \cite\NRSSD, Example 4.4, pp.\ 14--15\endcite\ for more details.   Looking ahead, in all these cases, Theorem \EESTARthm(b) remains true.   
\defRemark \PRODrem
\medbreak
\noindent
{\bf Remark \PRODrem.}\enspace Let $\big(B_1,\xbra\cdot\cdot_1,\|\cdot\|_1\big)$ be a Banach SSD space with Banach SSD dual\break $\big(B_1^*,\ybra\cdot\cdot_1,\|\cdot\|_1\big)$ and $\big(B_2,\xbra\cdot\cdot_2,\|\cdot\|_2\big)$ be a Banach SSD space with Banach SSD dual\break $\big(B_2^*,\ybra\cdot\cdot_2,\|\cdot\|_2\big)$.   We define $\|\cdot\|\colon\ B_1 \times B_2 \to \r$ and $\xbra\cdot\cdot\colon\ (B_1 \times B_2) \times (B_1 \times B_2) \to \r$ by $\big\|(b_1,b_2)\big\| := \sqrt{\|b_1\|_1^2 + \|b_2\|_2^2}$ and $\Xbra{(b_1,b_2)}{(c_1,c_2)} := \xbra{b_1}{c_1}_1 + \xbra{b_2}{c_2}_2$.   Similarly, we define $\|\cdot\|\colon\ B_1^* \times B_2^* \to \r$ and $\ybra\cdot\cdot\colon\ (B_1^* \times B_2^*) \times (B_1^* \times B_2^*) \to \r$ by $\big\|(b_1^*,b_2^*)\big\| := \sqrt{\|b_1^*\|_1^2 + \|b_2^*\|_2^2}$ and $\Ybra{(b_1^*,b_2^*)}{(c_1^*,c_2^*)} := \ybra{b_1^*}{c_1^*}_1 + \ybra{b_2^*}{c_2^*}_2$.   Then $\big(B_1 \times B_2,\xbra\cdot\cdot,\|\cdot\|\big)$ is a Banach SSD space with Banach SSD dual $\big(B_1^*\times B_2^*,\ybra\cdot\cdot,\|\cdot\|\big)$.
\smallskip
As an example of this construction, we could take  $\big(B_1,\xbra\cdot\cdot_1,\|\cdot\|_1\big)$ to be a Banach SSD space of the kind considered in Remark \Hrem, and $\big(B_2,\xbra\cdot\cdot_2,\|\cdot\|_2\big)$ to be a Banach SSD space of the kind considered in Example \EEDUALex.   If $B_1$ is odd-dimensional and $E$ is finite-dimensional then $B$ is odd--dimensional, and so cannot itself be of the form considered in Example \EEDUALex.   Example \Hex(c) is of this form, and a glance at that example shows how pathological the $q$--positive sets can be.   
\medbreak
We now recall Rockafellar's formula for the conjugate of a sum:
\defLemma \RSUMlem
\medbreak
\noindent
{\bf Lemma \RSUMlem.}\enspace\slant Let $X$ be a nonzero real Banach space and $f \in \PC(X)$, and let $h \in \PC(X)$ be real--valued and continuous.   Then, for all $x^* \in X^*$,
$$(f + h)^*(x^*) = \minn_{y^* \in X^*}\big[f^*(y^*) + h^*(x^* - y^*)\big].$$\par
\endslant
\Proof  See Rockafellar, \cite\FENCHEL, Theorem 3(a), p.\ 85\endcite, Z\u alinescu, \cite\ZBOOK, Theorem 2.8.7(iii), p.\ 127\endcite, or \cite\HBM, Corollary 10.3, p.\ 52\endcite.\qed
\defRemark \SHARPrem
\medbreak
\noindent
{\bf Remark \SHARPrem.}\enspace \cite\HBM, Theorem 7.4, p.\ 43\endcite\ contains a version of the Fenchel duality theorem with a sharp lower bound on the functional obtained.
\medbreak
\defLemma \SSDDlem
\medbreak
\noindent
{\bf Lemma \SSDDlem.}\enspace\slant Let $\big(B,\xbra\cdot\cdot,\|\cdot\|\big)$ be a Banach SSD space with Banach SSD dual\break $\big(B^*,\ybra\cdot\cdot,\|\cdot\|\big)$.   Define the function $g$ on $B$ by\quad $g := \half\|\cdot\|^2$.\quad   Let $f \in \PC(B)$.   Then
$$-\big((f - q) \episum p\big) = \big((f^* - \qt) \episum \pt\big)\circ \iota\ \on\ B.$$\endslant
\Proof Let $c \in B$.   Define $h\colon\ B \to \r$ by $h(b) := g(c - b)$.   Then, by direct computation using the fact that $g$ is an even function,
$$\all\ c^* \in B^*,\qquad h^*(c^*) = g^*(c^*) + \bra{c}{c^*}.\meqno\Hone$$
Then, from (\DUALthree), the continuity of $h$, Lemma \RSUMlem, (\Hone), (\DUALtwo) and the fact that, for all $c^* \in B^*$, $g^*(c^*) = \half\|c^*\|^2$,
$$\eqalign{-\big((f - q) \episum p\big)(c)
&= \supn_{b \in B}\big[-(f - q)(b) - p(c - b)\big]\cr
&= \supn_{b \in B}\big[\bra{b}{\iota(c)} -f(b) - h(b)\big] - q(c) = (f + h)^*\big(\iota(c)\big) - q(c)\cr
&=\minn_{b^* \in B^*}\big[f^*(b^*) + h^* \big(\iota(c) - b^*\big)\big] - q(c)\cr
&= \minn_{b^* \in B^*}\big[f^*(b^*) + g^*\big(\iota(c) - b^*\big) + \Bra{c}{\iota(c) - b^*}\big] - q(c)\cr
&= \minn_{b^* \in B^*}\big[f^*(b^*) + g^*\big(\iota(c) - b^*\big) - \ybra{\iota(c)}{b^*} + \qt\big(\iota(c)\big)\big]\cr
&= \minn_{b^* \in B^*}\big[(f^* - \qt)(b^*) + \pt\big(\iota(c) - b^*\big)\big]\cr
&= \big((f^* - \qt) \episum \pt\big)\big(\iota(c)\big).}$$
This completes the proof of Lemma \SSDDlem.\qed    
\defDefinition \MASdef
\medbreak
\noindent
{\bf Definition \MASdef.}\enspace Let $\big(B,\xbra\cdot\cdot,\|\cdot\|\big)$ be a Banach SSD space with Banach SSD dual\break $\big(B^*,\ybra\cdot\cdot,\|\cdot\|\big)$ and $f \in \PC(B)$.   We say that $f$ is an \slant MAS function\endslant\ if \quad $f \ge q$ on $B$\quad and \quad $f^* \ge \qt$ on $B^*$.\quad   This is an extension to Banach SSD spaces of the concept introduced by Marques Alves and Svaiter in \cite\ASTWO, Theorem 4.2, pp.\ 702--704\endcite\ for the situation described in Example \EEDUALex.
\medbreak
It is clear from the layout of Section \BSSDsec\ that the main results on VZ functions (that is, up to and including Theorem \EXthm) do not depend explicitly on $B^*$.   By contrast, a knowledge of $B^*$ is absolutely essential for even the definition of MAS function.   As a consequence, Theorem \MASVZthm\ below is rather suprising.   Theorem \MASVZthm(a) and its partial converse Theorem \MASVZthm(b) are motivated by various results scattered through \cite\VZ, Section 2\endcite.   Theorem \MASVZthm(c) is motivated by \cite\ASTHREE\endcite.
We recall from (\EEfive) that the $\pt$--density condition is satisfied in the situation of Example \EEDUALex, and we will discuss the implications of Theorem \MASVZthm\ to this example (via Theorem \EQthm) in Theorems \EEEQthm,  \SRthm, \BRthm\ and \DFPthm. 
\defTheorem \MASVZthm
\medbreak
\noindent
{\bf Theorem \MASVZthm.}\enspace\slant Let $\big(B,\xbra\cdot\cdot,\|\cdot\|\big)$ be a Banach SSD space with Banach SSD dual\break $\big(B^*,\ybra\cdot\cdot,\|\cdot\|\big)$.
\smallbreak
\noindent
{\rm(a)}\enspace Let $f \in \PC(B)$ be an MAS function.   Then $f$ is a VZ function.
\smallbreak
\noindent
{\rm(b)}\enspace
Let $\iota(B)$ be $\pt$--dense in $B^*$ and $f \in \PC(B)$ be a VZ function.   Then $f$ is an MAS function.
\smallbreak
\noindent
{\rm(c)}\enspace
Let $\iota(B)$ be $\pt$--dense in $B^*$ and $A \subset B$.   Then there exists an MAS function $f \in \PCLSC(B)$ such that $A = \PQ(f)\ \iff\ A$ is maximally $q$--positive and\quad $\Theta_A \ge \qt$ on $B^*$.\endslant
\Proof(a)\enspace Taking together (\IOTAfour), (\BSSDfour) and our hypothesis that $f$ is an MAS function, we have $\inf_B\big[f - q \big]\ge 0$, $\inf_Bp \ge 0$, $\inf_{B^*}\big[f^* - \qt\big] \ge 0$ and $\inf_{B^*}\pt \ge 0$.   Consequently, $\inf_B\big[(f - q) \episum p\big] \ge 0$ and $\inf_B\big[\big((f^* - \qt) \episum \pt\big)\circ \iota\big] \ge 0$, and (a) follows from Lemma \SSDDlem.
\smallbreak
(b)\enspace We know from (\VZthree) that\quad $f \ge q$ on $B$.\quad   Now let $b^* \in B^*$ and $c \in B$.   Then, from Lemma \SSDDlem\ again,
$$(f^* - \qt)(b^*) + \pt\big(\iota(c) - b^*\big)
\ge \big((f^* - \qt) \episum \pt\big)\big(\iota(c)\big)
= -\big((f - q) \episum p\big)(c) = 0.$$
Taking the infimum over $c \in B$ and using (\PSEUDOone),\quad $(f^* - \qt)(b^*) \ge 0$ on $B^*$.\quad   Since this holds for all $b^* \in B^*$, $f$ is an MAS function, giving (b).
\smallbreak 
(c)\enspace ($\lr$) Let $f \in \PCLSC(B)$ be an MAS function and $A = \PQ(f)$.   From (a), $f$ is a VZ function, and so Theorem \VZEXthm(d) implies that $A$ is maximally $q$--positive.   From (\VZONLYthree) and (b), $\Psi_A$ is an MAS function, consequently\quad ${\Psi_A}^* \ge \qt$ on $B^*$,\quad thus Lemma \THPSlem(c) implies that\quad $\Theta_A \ge \qt$ on $B^*$,\quad as required.
\smallskip
($\rl$) Suppose, conversely, that $A$ is maximally $q$--positive and\quad $\Theta_A \ge \qt$ on $B^*$.\quad From (\QPOSfive), Lemma \THPSlem(c) and (\POSFone),\quad $\Phi_A \ge q$ on $B$,\quad ${\Phi_A}^* \ge\qt$ on $B^*$\quad and $\PQ(\Phi_A) = A$, and the result follows with $f := \Phi_A$.   \big(We could also use $\Psi_A$ for this part of (c).\big)\qed       
\defDefinition \USEdef
\medbreak
\noindent
{\bf Definition \USEdef.}\enspace Let $\big(B,\xbra\cdot\cdot,\|\cdot\|\big)$ be a Banach SSD space with Banach SSD dual\break $\big(B^*,\ybra\cdot\cdot,\|\cdot\|\big)$.   We say that a topology $\T$ on $B^*$ is \slant compatible\endslant\ if it satisfies the conditions (a)--(c) below:
\smallbreak\noindent 
{\rm(a)}\enspace $\T \supset w(B^*,B^*)$.   {\rm ($w(B^*,B^*)$ is the weak topology induced on $B^*$ by the bilinear form $\ybra\cdot\cdot$.)} 
\smallbreak
\noindent
{\rm(b)}\enspace If $f \in \PCLSC(B)$ and $b^* \in B^*$ then there exists a net $\{b_\gamma\}$ of elements of $B$ such that\quad $\iota(b_\gamma) \to b^*$ in $\T$\quad and\quad $f(b_\gamma) \to f^{*@}(b^*)$.
\smallbreak
\noindent
{\rm(c)}\enspace If $\{b_\gamma\}$ and $\{a_\gamma\}$ are nets of elements of $B$, $b^* \in B^*$,\enspace $\iota(b_\gamma) \to b^*$ in $\T$\enspace and\enspace $\|a_\gamma - b_\gamma\| \to 0$ then\quad $\iota(a_\gamma) \to b^*$ in $\T$.
\defRemark \USErem
\medbreak
\noindent
{\bf Remark \USErem.}\enspace Definition \USEdef(a) says that $\T$ is not too small, Definition \USEdef(b) says that $\T$ is not too large, and Definition \USEdef(c) says that $\T$ behaves well under norm perturbations in $B$.   It follows from Lemma \TWOTOPlem\ below that $w(B^*,B^*)$ is compatible.   If the norm topology of $B^*$ is compatible then, from (b) above with $f := 0$, $\iota(B)$ is norm--dense in $B^*$.
\medbreak
Theorem \MASVZthm\ leads to the following fundamental result on the Gossez extension of a maximally $q$--positive set (see Definition \GOSSEZdef), which will be used in Theorem \EEEQthm, and thus indirectly in Theorems \BRthm\ and \DFPthm, which depend on Theorem \EEEQthm.   It is actually Theorem \EQthm\ that provides the incentive for the investigation of the continuity of $\qt$ that we will perform in Section \QTsec.
\defTheorem \EQthm
\medbreak
\noindent
{\bf Theorem \EQthm}\enspace\slant Let $\big(B,\xbra\cdot\cdot,\|\cdot\|\big)$ be a Banach SSD space with Banach SSD dual\break $\big(B^*,\ybra\cdot\cdot,\|\cdot\|\big)$, $\T$ be a compatible topology on $B^*$, $\qt$ be $\T$--continuous and $A$ be a maximally $q$--positive subset of $B$.  Then the conditions {\rm(a)--(c)} below are equivalent.
\smallbreak\noindent
{\rm(a)}\enspace $\AGOSS \subset \iota(A)^\T$, the closure of $\iota(A)$ in the topology $\T$.\smallbreak\noindent
{\rm(b)}\enspace For all $b^* \in \AGOSS$, $\inf\qt\big(\iota(A) - b^*\big) \le 0$.
\smallbreak\noindent
{\rm(c)}\enspace $\Theta_A \ge \qt$ on $B^*$.\endslant
\Proof Suppose that (a) is satisfied and $b^* \in \AGOSS$.   Then there exists a net $\{a_\gamma\}$ of elements of $A$ such that\quad $\iota(a_\gamma) \to b^*$ in $\T$.\quad From Definition \USEdef(a), $\iota(a_\gamma) \to b^*$ in $w(B^*,B^*)$ and so $\Ybra{\iota(a_\gamma)}{b^*} \to \Ybra{b^*}{b^*} = 2\qt(b^*)$.
From the $\T$--continuity of $\qt$,\quad $\qt\big(\iota(a_\gamma)\big) \to \qt(b^*)$.\quad Thus
$$\qt\big(\iota(a_\gamma) - b^*\big) = \qt\big(\iota(a_\gamma)\big) - \Ybra{\iota(a_\gamma)}{b^*} + \qt(b^*) \to \qt(b^*) - 2\qt(b^*) + \qt(b^*) = 0,$$ 
and so (a)$\ \lr\ $(b).   It follows from Theorem \AGOSSthm(b) that (b)$\ \lr\ $(c).   So it remains to prove that (c)$\ \lr\ $(a).
\smallbreak
So suppose that (c) is satisfied and $b^* \in \AGOSS$.   (\QPOSfour) implies that $\Phi_A \in \PCLSC(B)$.   From Theorem \AGOSSthm(c), ${\Phi_A}^{*@}(b^*) = \qt(b^*)$.   Definition \USEdef(b) now gives us a net $\{b_\gamma\}$ of elements of $B$ such that
$$\iota(b_\gamma) \to b^*\ \hbox{in}\ \T \quand \Phi_A(b_\gamma) \to {\Phi_A}^{*@}(b^*) = \qt(b^*).\meqno\LIMthree$$
It now follows from (\DUALtwo), the first assertion in (\LIMthree) and the $\T$--continuity of $\qt$ that $q(b_\gamma) = \qt\circ\iota(b_\gamma) \to \qt(b^*)$ and so, using the second assertion in (\LIMthree),
$$(\Phi_A - q)(b_\gamma) = \Phi_A(b_\gamma) - q(b_\gamma) \to \qt(b^*) - \qt(b^*) = 0.\meqno\LIMfour$$
Now (\QPOSfive) implies that $\Phi_A \ge q$ on $B$,\quad and, from Lemma \THPSlem(c),\quad ${\Phi_A}^* \ge \Theta_A \ge\qt$ on $B^*$, from which $\Phi_A$ is an MAS fucntion.\quad Thus, from Theorem \MASVZthm(a), $\Phi_A$ is a VZ function.   Since $\Phi_A$ is lower semicontinuous on $B$, (\POSFone) and (\FPHItwo) imply that, for all $\gamma$,
$$\dist(b_\gamma,A) = \dist\big(b_\gamma,\PQ(\Phi_A)\big) \le \rttwo\sqrt{(\Phi_A - q)(b_\gamma)}.$$
\big(We interpret $\sqrt\infty$ to be $\infty$.\big)   (\LIMfour) now gives us $a_\gamma \in A$ such that $\|a_\gamma - b_\gamma\| \to 0$, and so (a) follows from the first assertion in (\LIMthree) and Definition \USEdef(c).\qed
\defSection \CLBsec
\medbreak
\leftline{\bf \CLBsec.\quad $\CLB(X)$ and $\TCLB(X\dbs)$}
\medskip
\noindent
Let $X$ be a nonzero real Banach space.   Corresponding to the usage outlined in first paragraph of Section \EXTsec, if $f \in \PC(X)$ and $f^* \in \PC(X^*)$, we define $f\dbs\colon X\dbs \to \rbar$ by $f\dbs(x\dbs) := \sup_{X^*}\big[\bra{\cdot}{x\dbs} - f^*\big]$.   We write $\CLB(X)$ for the set of all convex functions $f\colon\ X \to \r$ that are Lipschitz on the bounded subsets of $X$, or equivalently, bounded above on the bounded subsets of $X$, and we define the topology $\TCLB(X\dbs)$ on $X\dbs$ to be the coarsest topology on $X\dbs$ making all the functions\quad$h\dbs\colon\ X\dbs \to \r \quad \big(h \in \CLB(X)\big)$\quad continuous.   \big(See \cite\HBM, Definition 38.1, p,\ 155\endcite.\big)   We write $\TNORM(X)$ for the norm--topology on $X$.
We collect together in the following lemma the basic properties of $\TCLB(X\dbs)$ that we will use.   We will discuss subtler properties of the topologies $\TCLB$ in Lemma \TCLBElem.
\defLemma \TCBlem
\medbreak
\noindent
{\bf Lemma \TCBlem.}\enspace\slant Let $X$ be a nonzero real Banach space.\smallbreak\noindent
{\rm(a)}\enspace Let $\big\{x_\gamma\dbs\big\}$ be a net of elements of $X\dbs$, $x\dbs \in X\dbs$ and $x_\gamma\dbs \to x\dbs$ in $\TCLB(X\dbs)$.  Then $\big\{x_\gamma\dbs\big\}$ is eventually bounded and $x_\gamma\dbs \to x\dbs$ in the weak$^*$--topology $w(X\dbs,X^*)$.\smallbreak
\noindent
{\rm(b)}\enspace Let $f \in \PCLSC(X)$  and $x\dbs \in X\dbs$.   Then there exists a net $\{x_\gamma\}$ of elements of $X$ such that\quad $\wh{x_\gamma} \to x\dbs$ in $\TCLB(X\dbs)$\quad and\quad $f(x_\gamma) \to f\dbs(x\dbs)$.   {\rm(We note from the Fenchel--Moreau theorem that $f^* \in \PC(X^*)$.)}
\smallbreak
\noindent
{\rm(c)}\enspace Let $\{x\dbs_\gamma\}$ and $\{y\dbs_\gamma\}$ be nets of elements of $X\dbs$, $x\dbs \in X\dbs$,\quad $x\dbs_\gamma \to x\dbs$ in $\TCLB(X\dbs)$\quad and $\|y\dbs_\gamma - x\dbs_\gamma\| \to 0$.\quad   Then\quad $y\dbs_\gamma \to x\dbs$ in $\TCLB(X\dbs)$.\endslant
\Proof(a)\enspace See \cite\HBM, Lemma 38.2(b,f), p.\ 156\endcite.
\smallbreak
(b)\enspace See \cite\HBM, Lemma 45.9(a), p.\ 175\endcite.
\smallbreak
(c)\enspace See \cite\HBM, Lemma 45.15, p. 177\endcite.\qed
\defRemark \PATHOrem
\medbreak
\noindent
{\bf Remark \PATHOrem.}\enspace Despite the nice properties of $\TCLB(X\dbs)$ exhibited in Lemma \TCBlem, it is nevertheless quite a pathological topology.   For instance, if $\big(B\dbs,\TCLB(X\dbs)\big)$ is a topological vector space then $X$ is reflexive.   See \cite\HBM, Remark 45.13, p. 177\endcite.
\medbreak  
Lemma \TCLBElem(b) was originally developed in a study of the subdifferentials of saddle\break functions.  
\defLemma \TCLBElem
\medbreak
\noindent
{\bf Lemma \TCLBElem.}\enspace\slant Let $E$ be a nonzero Banach space.
\smallbreak
\noindent
{\rm(a)}\enspace The map $\qt$ is continuous from $\big(E^* \times E\dbs,\TNORM(E^*) \times \TCLB(E\dbs)\big)$ into $\r$.
\smallbreak
\noindent
{\rm(b)}\enspace Let $H$ also be a nonzero Banach space, $(y\dbs,z) \in E\dbs \times H$ and $\big\{(y_\gamma\dbs,z_\gamma)\big\}$ be a net of elements of $E\dbs \times H$.   Then\quad $(y_\gamma\dbs,\wh{z_\gamma}) \to (y\dbs,\wh{z})$ in $\TCLB(E\dbs \times H\dbs) \iff\break(y_\gamma\dbs,z_\gamma) \to (y\dbs,z)$ in $\TCLB(E\dbs) \times \TNORM(H)$.\endslant
\Proof(a)\enspace See \cite\HBM, Lemma 38.2(e), p.\ 156\endcite.   \big(We note that $B^*$ in that reference was defined to be $E\dbs \times E^*$ rather than $E^* \times E\dbs$ as we have done here, and the topology $\TCLBN(B^*)$ was defined to be $\TCLB(E\dbs) \times \TNORM(E^*)$.\big)
\smallbreak
(b)\enspace See \cite\HBM, Theorem 49.4, pp.\ 194\endcite.\qed
\medbreak
\defSection \QTsec
\medbreak
\leftline{\bf \QTsec.\quad $\TD(B^*)$}
\medskip
\noindent
We suppose throughout this section that $\big(B,\xbra\cdot\cdot,\|\cdot\|\big)$ is a Banach SSD space with Banach SSD dual $\big(B^*,\ybra\cdot\cdot,\|\cdot\|\big)$.   In order to apply Theorem \EQthm, we need a compatible topology on $B^*$ with respect to which $\qt$ is continuous.      To get some insight into this problem, we consider the case of Example \EEDUALex, that is to say, $B^* = E^* \times E\dbs$ and $\qt\colon (x^*,x\dbs)\mapsto \bra{x^*}{x\dbs}$.   It has been known since Gossez's work in \cite\GOSSEZ\endcite\ that $\TNORM(E^* \times E\dbs)$ is too large to be of any practical use.   (The root of the problem can be found in Remark \USErem.)   Gossez considers the topology $\TNORM(E^*) \times w(E\dbs,E^*)$, but this topology does not seem to generalize easily to the case of SSD spaces.   In Definition \TDdef, we introduce the topology $\TD(B^*)$ on $B^*$.   We will see in Lemma \TWOTOPlem\ that $\TD(B^*)$ is sufficiently small that it is compatible, and we will see in Theorem \EESTARthm(b) that $\TD(B^*)$ is sufficiently large that Theorem \EQthm\ leads to significant results.
\defDefinition \TDdef
\medbreak
\noindent
{\bf Definition \TDdef.}\enspace Let $\big(B,\xbra\cdot\cdot,\|\cdot\|\big)$ be a Banach SSD space with Banach SSD dual $\big(B^*,\ybra\cdot\cdot,\|\cdot\|\big)$.   We define the topology $\TD(B^*)$ on $B^*$ to be the coarsest topology on $B^*$ making the function\quad$\it\colon\ B^* \to \big(B\dbs,\TCLB(B\dbs)\big)$\quad continuous.   This means that if $\{b^*_\gamma\}$ is a net of elements of $B^*$ and $b^* \in B^*$   then
$$b^*_\gamma \to b^*\ \hbox{in}\ \TD(B^*) \iff \it(b^*_\gamma) \to \it(b^*)\ \hbox{in}\ \TCLB(B\dbs).\meqno\TDone$$
Now suppose that $\{b_\gamma\}$ is a net of elements of $B$ and $b^* \in B^*$.   Combining (\TDone) with Lemma \FACTORlem, we have
$$\iota(b_\gamma) \to b^*\ \hbox{in}\ \TD(B^*) \iff \wh{b_\gamma} \to \it(b^*)\ \hbox{in}\ \TCLB(B\dbs).\meqno\TDtwo$$\par    
\defLemma \TWOTOPlem
\medbreak
\noindent
{\bf Lemma \TWOTOPlem.}\enspace\slant We suppose that $\big(B,\xbra\cdot\cdot,\|\cdot\|\big)$ is a Banach SSD space with Banach SSD dual $\big(B^*,\ybra\cdot\cdot,\|\cdot\|\big)$.   Then $\TD(B^*)$ is a compatible topology on $B^*$.\endslant
\Proof We first verify Definition \USEdef(a).   Let $\{b_\gamma^*\}$ be a net of elements of $B^*$, $b^* \in B^*$ and\quad $b_\gamma^* \to b^*$ in $\TD(B^*)$.\quad  (\TDone) implies that\quad $\it(b^*_\gamma) \to \it(b^*)$ in $\TCLB(B\dbs)$\quad and so, from Lemma \TCBlem(a),\quad $\it(b^*_\gamma) \to \it(b^*)$ in $w(B\dbs,B^*)$.\quad (\BSSDthree) now gives us that
$$\all\ c^* \in B^*,\quad\ybra{c^*}{b^*_\gamma} = \Bra{c^*}{\it(b^*_\gamma)} \to \Bra{c^*}{\it(b^*)} = \ybra{c^*}{b^*},$$
and so $b_\gamma^* \to b^*$ in $w(B^*,B^*)$.   This completes the proof of Definition \USEdef(a). 
\smallbreak
We next verify Definition \USEdef(b).   To this end, let $f \in \PCLSC(B)$ and $b^* \in B^*$.  Lemma \TCBlem(b) provides us with a net $\{b_\gamma\}$ of elements of $B$ such that\quad $\wh{b_\gamma} \to \it(b^*)$ in $\TCLB(B\dbs)$\quad and\quad $f(b_\gamma) \to f\dbs\circ\it(b^*)$.\quad (\TDtwo) gives\quad $\iota(b_\gamma) \to b^*$ in $\TD(B^*)$,\quad  and the analog of (\LINKthree) gives \quad $f\dbs\circ\it(b^*) = f^{*@}(b^*)$.\quad This completes the proof of Definition \USEdef(b).
\smallbreak
Finally, we verify Definition \USEdef(c).   To this end, let $\{b_\gamma\}$ and $\{a_\gamma\}$ be nets of elements of $B$, $b^* \in B^*$,\quad $\iota(b_\gamma) \to b^*$ in $\TD(B^*)$\quad and\quad $\|a_\gamma - b_\gamma\| \to 0$.\quad   From  (\TDtwo),\quad $\wh{b_\gamma} \to \it(b^*)$ in $\TCLB(B\dbs)$.\quad   Since $\ \wh{}\ $ is a norm--isometry, $\|\wh{a_\gamma} - \wh{b_\gamma}\| \to 0$, and so Lemma \TCBlem(c) implies that \quad $\wh{a_\gamma} \to \it(b^*)$ in $\TCLB(B\dbs)$.\quad  It now follows from another application of (\TDtwo) that $\iota(a_\gamma) \to b^*$ in $\TD(B^*)$.   This completes the proof of Definition \USEdef(c). \big(In fact, one can prove in a similar way, using (\TDone) instead of  (\TDtwo), the stronger result that if $\{b_\gamma^*\}$ and $\{a_\gamma^*\}$ are nets of elements of $B^*$, $b^* \in B^*$,\quad $b_\gamma^* \to b^*$ in $\TD(B^*)$\quad and\quad $\|a_\gamma^* - b_\gamma^*\| \to 0$ then\quad $a_\gamma^* \to b^*$ in $\TD(B^*)$.\big)
\qed     
\medbreak
Theorem \EESTARthm\ below will be used in Theorems \EEEQthm, \BRthm, and \DFPthm.  
\defTheorem \EESTARthm
\medbreak
\noindent
{\bf Theorem \EESTARthm}\enspace\slant Let $E$ be a nonzero Banach space and $\big(B,\xbra\cdot\cdot,\|\cdot\|_2\big)$ and $\big(B^*,\ybra\cdot\cdot,\|\cdot\|_2\big)$ be as in Example \EEDUALex.
\smallbreak
\noindent
{\rm(a)}\enspace The topologies $\TD(B^*)$ and $\TNORM(E^*) \times \TCLB(E\dbs)$ on $B^* = E^* \times E\dbs$ are identical.  
\smallbreak
\noindent
{\rm(b)}\enspace $\qt$ is $\TD(B^*)$--continuous. 
\endslant
\Proof We recall from  Example \EEDUALex\ that, for all $(y^*,y\dbs) \in B^*$, $\it(y^*,y\dbs) = \big(y\dbs,\wh{y^*}\big)$ and $\qt(y^*,y\dbs) = \bra{y^*}{y\dbs}$.   Let $\big\{(y^*_\gamma,y_\gamma\dbs)\big\}$ be a net of elements of $B^*$ and $(y^*,y\dbs) \in B^*$.   Then, from (\TDone),
$$\eqalign{(y^*_\gamma,y_\gamma\dbs) \to (y^*,y\dbs)\ \hbox{in}\ \TD(B^*)\iff
\big(y_\gamma\dbs,\wh{y^*_\gamma}\big) \to \big(y\dbs,\wh{y^*}\big)\ \hbox{in}\ \TCLB(B\dbs).}$$
(a) is now immediate from Lemma \TCLBElem(b) with $H := E^*$, and (b) is immediate from (a) and Lemma \TCLBElem(a).\qed     
\defRemark \BONUSrem
\medbreak
\noindent
{\bf Remark \BONUSrem.}\enspace A hidden bonus of Theorem \EESTARthm\ is that, despite the fact that $B\dbs =\break E\dbs \times E^{***}$, we do not actually have to deal with $E^{***}$.
\defSection \MONsec
\medbreak
\leftline{\bf \MONsec.\quad Classes of monotone sets}
\medskip
\noindent
We suppose in this section that $E$ is a nonzero Banach space.   For most of the time it will be convenient to work in terms of subsets of $E \times E^*$ and $E^* \times E\dbs$ rather than multifunctions $E \toto E^*$ and $E\dbs \toto E^*$, and we leave it to the reader to verify the consistence between our versions and the multifunction versions. 
\medbreak
The motivation for the consideration of the various classes of sets described below was to see how many of the properties of maximally monotone sets on reflexive spaces can be recovered in the nonreflexive case.   Historically, the first such classes of sets were the class of sets of ``dense type'' and ``type (D)''.   These were essentially introduced by Gossez in \cite\GOSSEZ, Lemme 2.1, p.\ 375\endcite\ --- see Phelps, \cite\PRAGUE, Section 3\endcite\ for an exposition.   If $A$ is a monotone subset of $E \times E^*$, Gossez defines $\Abar \subset\ E^* \times E\dbs$ by:
$$\Abar := \big\{(y^*,y\dbs) \in E^* \times E\dbs\colon\ \infn_{(a,a^*) \in A}\bra{a^* - y^*}{\wh{a} - y\dbs} \ge 0\big\}.\meqno\ABARone$$\par
\defDefinition \TYPEDdef
\medbreak
\noindent
{\bf Definition \TYPEDdef.}\enspace Let $A \subset E \times E^*$.   We say that $A$ is \slant maximally monotone of type (D)\endslant\ if $A$ is maximally monotone and, for all $(y^*,y\dbs) \in \Abar$, there exists a bounded net $\{(a_\gamma,a^*_\gamma)\}$ of elements of $A$ such that  $(a^*_\gamma,\wh{a_\gamma}) \to (y^*,y\dbs)$ in $\TNORM(E^*) \times w(E\dbs,E^*)$.   We say that $A$ is \slant maximally monotone of dense type\endslant\ if the topology $w(E\dbs,E^*)$ in the definition above is replaced by the topology ${\cal T}_1$ defined to be the upper bound of $w(E\dbs,E^*)$ and the coarsest topology making the function $\|\cdot\|\colon\ E\dbs \to \r$ continuous.    
\medbreak
The next classes of monotone sets in our discussion are the classes of sets of type (NI) and (WD), which were introduced in \cite\RANGE, Definition 10, p.\ 183\endcite\ and \cite\RANGE, Definition 14, p.\ 187\endcite.
\defDefinition \NIdef
\medbreak
\noindent
{\bf Definition \NIdef.}\enspace Let $A \subset E \times E^*$.   We say that $A$ is \slant maximally monotone of type (NI)\endslant\ if $A$ is maximally monotone and,
$$\all\ (y^*,y\dbs) \in E^* \times E\dbs,\quad \infn_{(a,a^*) \in A}\bra{a^* - y^*}{\wh a - y\dbs} \le 0.$$
We say that $A$ is \slant maximally monotone of type (WD)\endslant\ if $A$ is maximally monotone and, for all $(y^*,y\dbs) \in \Abar$, there exists a bounded net $\{(a_\gamma,a^*_\gamma)\}$ of elements of $A$ such that $a^*_\gamma \to y^*$ in $\TNORM(E^*)$.   Clearly, 
$$\slant\hbox{if}\ A\ \hbox{is maximally monotone of type (D) then}\ A\ \hbox{is of type (WD)},\endslant\meqno\NIone$$
and it was proved in \cite\RANGE, Lemma 15, pp.\ 187--188\endcite\ that
$$\slant\hbox{if}\ A\ \hbox{is maximally monotone of type (WD) then}\ A\ \hbox{is of type (NI)}\endslant.\meqno\NItwo$$
\par
The next class of monotone sets in our discussion is the class of sets of type (ED), which was introduced in \cite\MANDM, Definition 35.1, p.\ 138\endcite\ under the name ``type (DS)''.
\defDefinition \TYPEEDdef
\medbreak
\noindent
{\bf Definition \TYPEEDdef.}\enspace Let $A \subset E \times E^*$.    We say that $A$ is \slant maximally monotone of type (ED)\endslant\ if $A$ is maximally monotone and, for all $(y^*,y\dbs) \in \Abar$, there exists a net $\{(a_\gamma,a^*_\gamma)\}$ of elements of $A$ such that
$(a^*_\gamma,\wh{a_\gamma}) \to (y^*,y\dbs)$ in $\TNORM(E^*) \times \TCLB(E\dbs)$.   It is clear from Lemma \TCBlem(a) and Definition \TYPEDdef\ that
$$\slant\hbox{if}\ A\ \hbox{is maximally monotone of type (ED) then}\ A\ \hbox{is of dense type}\endslant.\meqno\NIthree$$
\par
We now recast the above definitions in the more compact notation of SSD spaces, using the notation of Example \EEDUALex.
\defLemma \RECASTlem
\medbreak
\noindent
{\bf Lemma \RECASTlem}\enspace\slant Let $A$ be a maximally monotone subset of $E \times E^*$.   Then:\par
\noindent
{\rm(a)}\enspace $\Abar = \AGOSS$.
\smallbreak
\noindent
{\rm(b)}\enspace $A$ is of type (NI) $\iff \Theta_A \ge \qt$ on $E^* \times E\dbs$.
\smallbreak
\noindent
{\rm(c)}\enspace $A$ is of type (ED) $\iff$ for all $(y^*,y\dbs) \in \AGOSS$, there exists a net $\{(a_\gamma,a^*_\gamma)\}$ of elements of $A$ such that\quad $\iota(a_\gamma,a^*_\gamma) \to (y^*,y\dbs)$ in $\TD(E^* \times E\dbs)$. 
\endslant
\Proof(a) is immediate from (\ABARone), (\EEeight) and the first expression in (\GOSSEZone), and (b) is immediate from Definition \NIdef, (\EEeight) and Lemma \HOMlem(b).   As for (c), from Definition \TYPEEDdef\ and (a), $A$ is of type (ED) exactly when, for all $(y^*,y\dbs) \in \AGOSS$, there exists a net $\{(a_\gamma,a^*_\gamma)\}$ of elements of $A$ such that
$(a^*_\gamma,\wh{a_\gamma}) \to (y^*,y\dbs)$ in $\TNORM(E^*) \times \TCLB(E\dbs)$ and (c) follows from Theorem \EESTARthm(a).\qed   
\medbreak
Now it is clear from  (\NIthree), (\NIone) and (\NItwo) that, for maximally monotone sets,
$$\slant\hbox{type (ED)}\qlr\hbox{dense type}\qlr\hbox{type (D)}\qlr\hbox{type (WD)}\qlr\hbox{type (NI)},\endslant$$
and the question arises naturally whether there are any result in the reverse direction.   Considerable progress was made recently by Marques Alves and Svaiter in \cite\ASD, Theorem 4.4, pp.\ 10--11\endcite, where it was established that    
$$\slant\hbox{if}\ A\ \hbox{is maximally monotone of type (NI) then}\ A\ \hbox{is of type (D)},\endslant\meqno\SRone$$
thus for maximally monotone sets, type (D), type (WD) and type (NI) are equivalent.   Consequently, the conjecture on \cite\RANGE, p.\ 187\endcite\ and the first conjecture on \cite\RANGE, p.\ 188\endcite\ are false while, as we will see in Theorem \BRthm(f), the second conjecture on \cite\RANGE, p.\ 188\endcite\ (on the convexity of the closure of the range) is true.   This latter result was actually established by Zagrodny in \cite\ZAGRODNY\endcite\ \big(see \cite\HBM, Problem 43.3, p.\ 168\endcite\big).   (\SRone) also provides a positive answer to \cite \HBM, Problem 36.4, p.\ 149\endcite.   The following result extends (\SRone), and provides an (unexpected positive) answer to \cite\BR, Problem 4.3, p.\ 268\endcite:
\defTheorem \EEEQthm
\medbreak
\noindent
{\bf Theorem \EEEQthm}\enspace\slant Let $E$ be a nonzero Banach space. Then for maximally monotone subsets of $E \times E^*$, type (ED), dense type, type (D), type (WD) and type (NI) are equivalent.\endslant
\Proof By virtue of the remarks above, we only have to prove that\quad type (NI)$\qlr$type (ED).\quad   So let $A$ be a maximally monotone subset of $E \times E^*$ of type (NI).  Lemma \RECASTlem(b) implies that $\Theta_A \ge \qt$ on $E^* \times E\dbs$, and then, from Lemma \RECASTlem(a) and Theorem \EQthm\big((c)$\lr$(a)\big), for all $(y^*,y\dbs) \in \AGOSS$, there exists a net $\{(a_\gamma,a^*_\gamma)\}$ of elements of $A$ such that\quad $\iota(a_\gamma,a^*_\gamma) \to (y^*,y\dbs)$ in $\TD(E^* \times E\dbs)$.    Thus, from Lemma \RECASTlem(c), $A$ is of type\break (ED).\qed
\medbreak
The next class of monotone sets in our discussion is the class of strongly representable sets, which was introduced and studied in \cite\ASTWO\endcite, \cite\ASTHREE\endcite\ and \cite\VZ\endcite. 
\defDefinition \SRdef
\medbreak
\noindent
{\bf Definition \SRdef.}\enspace Let $E$ be a nonzero Banach space and $A \subset E \times E^*$.   We say that $A$ is \slant strongly representable\endslant\ if there exists an MAS function $f \in \PCLSC(E \times E^*)$ such that $A = \PQ(f)$.
\medbreak
We now give a proof using SSD spaces of the following result, which was established by Marques Alves and Svaiter in \cite\ASTWO, Theorem 4.2, pp.\ 702--704\endcite\ and \cite\ASTHREE, Theorem 1.2\endcite.
\defTheorem \SRthm
\medbreak
\noindent
{\bf Theorem \SRthm.}\enspace\slant Let $E$ be a nonzero Banach space and $A \subset E \times E^*$.   Then $A$ is strongly representable $\iff A$\ is maximally monotone of type (NI).\endslant
\Proof This is immediate from (\EEfive), Theorem \MASVZthm(c) and Lemma \RECASTlem(b).\qed
\medbreak
The ``maximally monotone'' assertion of Theorem \BRthm(a) was obtained in \cite\VZ, Theorem 8\endcite\ under the VZ hypothesis and, in \cite\ASTWO, Theorem 4.2(2)\endcite\ under the MAS hypothesis.
\smallskip
Theorem \BRthm(c) extends the result proved in \cite\VZ, Corollary 25\endcite\ that $\PQ(f)$ is of type (ANA).
\smallskip
Theorem \BRthm(d) extends the result proved in \cite\ASTWO, Theorem 4.2(2)\endcite.
\smallskip
Theorem \BRthm(f) was obtained in \cite\VZ, Corollary 7\endcite.   This is a very significant result, because maximally monotone sets $A$ of $E \times E^*$ are known such that $\overline{\pi_{E^*}(A)}$ is not convex.   \big(The first such example was given by Gossez in \cite\GOSSEZC, Proposition, p. 360\endcite\big).   Thus \big(as was first observed in \cite\VZ\endcite\big) Theorem \BRthm(f) implies that there exist maximally monotone sets $A$ that are not of the form $\PQ(f)$ for any lower semicontinuous VZ function on $E \times E^*$ or, equivalently, not of the form $\PQ(f)$ for any lower semicontinuous MAS function on $E \times E^*$.
\smallskip
In \cite\ZAGRODNY, Section 3, pp.\ 775--783\endcite, Zagrodny considers subsets $S$ of $E \times E^*$ such that, writing $\phi_S$ for the Fitzpatrick function of $S$, $\phi_S \in \PC(E \times E^*)$,\quad $\phi_S \ge q$ on $E \times E^*$\quad and
$$(a^*,a\dbs) \in \partial\phi_S(x,x^*) \qlr \phi_S(x,x^*) \le \bra{x}{a^*} + \bra{x^*}{a\dbs} - \bra{a^*}{a\dbs}.\meqno\ZAGone$$
Since the analysis in \cite\ZAGRODNY\endcite\ leans heavily on $\eps$--enlargements, it is hard to correlate it on a step by step basis with what we have presented here.   Nevertheless, we note the following consequences if $A$ is a maximally monotone subset of $E \times E^*$ of type (NI):
\par
\noindent
$\bullet$\enspace In \cite\ZAGRODNY, (20), p.\ 776\endcite, Zagrodny deduces the second assertion in (\NAPAIRone).
\par
\noindent
$\bullet$\enspace In \cite\ZAGRODNY, Corollary 3.4, p.\ 780\endcite, Zagrodny deduces that $A$ is $p$--dense in $E \times E^*$ and an additional boundedness conclusion on the approximants.   \big(Compare Theorem \VZEXthm(b).\big)
\par
\noindent
$\bullet$\enspace In \cite\ZAGRODNY, Corollary 3.5, pp.\ 781--782\endcite\ and \cite\ZAGRODNY, Corollary 3.6, pp.\ 782--783\endcite, Zagrodny deduces that the sets $\overline{\pi_E(A)}$ and $\overline{\pi_{E^*}(A)}$ are convex.   As we have already noted, this is a very significant result.
\defProblem \ZAGprob
\medbreak
\noindent
{\bf Problem \ZAGprob.}\enspace We note from the definition of $\partial \phi_S$ that (\ZAGone) can be put in the form\quad $(a^*,a\dbs) \in \partial\phi_S(x,x^*) \lr {\phi_S}^*(a^*,a\dbs) \ge \bra{a^*}{a\dbs}$\quad that is to say\quad $(a^*,a\dbs) \in R(\partial\phi_S)\break\lr {\phi_S}^*(a^*,a\dbs) \ge \bra{a^*}{a\dbs}$.\quad   This leads to the following question:  \slant is Theorem \MASVZthm(a) true if, instead of assuming that $f$ is an MAS function, we assume that\quad $f \ge q$ on $B$\quad and\quad $f^* \ge \qt$ on $R(\partial f)$?\endslant\quad   Given the applications of Theorem \MASVZthm(a) that we make, it is probably no restriction to assume that $f \in \PCLSC(B)$.    
\defTheorem \BRthm
\medbreak
\noindent
{\bf Theorem \BRthm.}\enspace\slant Let $E$ be a nonzero Banach space and $f \in \PCLSC(E \times E^*)$.   Assume either that $f$ is a VZ function or, equivalently {\rm\big(bearing in mind (\EEfive) and  Theorem \MASVZthm(a,b)\big)}, an MAS function, and let $A := \PQ(f)$.   Then:
\par\noindent
{\rm(a)}\enspace $A$ is a maximally monotone subset of $E \times E^*$ of type (ED).
\smallskip\noindent
{\rm(b)}\enspace Let $(x,x^*) \in E \times E^*$ and $\alpha,\beta > 0$.   Then there exists a unique value of $\tau \ge 0$ for which there exists a bounded sequence $\big\{(y_n,y_n^*)\big\}_{n \ge 1}$ of elements of $A$ such that,
$$\lim_{n \to \infty}\|y_n - x\| = \alpha\tau,\quad \lim_{n \to \infty}\|y_n^* - x^*\| = \beta\tau \quand \lim_{n \to \infty}\bra{y_n - x}{y_n^* - x^*} = - \alpha\beta\tau^2.$$
{\rm(c)}\enspace Let $(x,x^*) \in E \times E^* \setminus A$ and $\alpha,\beta > 0$.   Then there exists a bounded sequence $\big\{(y_n,y_n^*)\big\}_{n \ge 1}$ of elements of\quad $A \cap \big[(E \setminus \{x\}) \times (E^* \setminus \{x^*\})\big]$\quad such that,
$$\lim_{n \to \infty}\f{\|y_n - x\|}{\|y_n^* - x^*\|} = \f\alpha\beta \quand \lim_{n \to \infty}\f{\bra{y_n - x}{y_n^* - x^*}}{\|y_n - x\|\|y_n^* - x^*\|} = -1.\meqno\NAPAIRone$$
In particular, $A$ is of type (ANA) {\rm\big(see \cite\HBM, Definition 36.11, p.\ 152\endcite\big)}.
\smallskip\noindent
{\rm(d)}\enspace Let $(x,x^*) \in E \times E^* \setminus A$, $\alpha,\beta > 0$ and $\inf_{(y,y^*) \in A}\bra{y - x}{y^* - x^*} > -\alpha\beta$.  Then there exists a bounded sequence $\big\{(y_n,y_n^*)\big\}_{n \ge 1}$ in\quad $A \cap \big[(E \setminus \{x\}) \times (E^* \setminus \{x^*\})\big]$\quad  such that {\rm(\NAPAIRone)} is satisfied, $\lim_{n \to \infty}\|y_n - x\| < \alpha$ and $\lim_{n \to \infty}\|y_n^* - x^*\| < \beta$.   In particular, $A$ is of type (BR) {\rm\big(see \cite\HBM, Definition 36.13, p.\ 153\endcite\big)}.
\smallskip\noindent
{\rm(e)}\enspace Let $(x,x^*) \in E \times E^* \setminus A$, $\alpha,\beta > 0$ and $f(x,x^*) < \bra{x}{x^*} + \alpha\beta$.   Then there exists a bounded sequence $\big\{(y_n,y_n^*)\big\}_{n \ge 1}$ of elements of \quad $A \cap \big[(E \setminus \{x\}) \times (E^* \setminus \{x^*\})\big]$\quad  such that {\rm(\NAPAIRone)} is satisfied, $\lim_{n \to \infty}\|y_n - x\| < \alpha$ and $\lim_{n \to \infty}\|y_n^* - x^*\| < \beta$.
\smallskip\noindent
{\rm(f)}\enspace We define the {\rm projection maps} $\pi_E\colon E \times E^* \to E$ and $\pi_{E^*}\colon E \times E^* \to E^*$ by $\pi_E(x,x^*) := x$ and $\pi_{E^*}(x,x^*) := x^*$.   Then\quad $\overline{\pi_E(A)} = \overline{\pi_E(\dom\,f)}$ and $\overline{\pi_{E^*}(A)} = \overline{\pi_{E^*}(\dom\,f)}$.   Consequently, the sets $\overline{\pi_E(A)}$ and $\overline{\pi_{E^*}(A)}$ are convex.
\endslant
\Proof(a) is immediate from Theorems \VZEXthm(d), \SRthm\ and \EEEQthm.
\smallbreak
(b), (c) and (d) are immediate from (a) and either \cite\BR, Theorem 8.6, pp.\ 277--278\endcite\ or \cite\HBM, Theorem 42.6, pp.\ 163--164\endcite.
\smallbreak
(e) is immediate from (d) and the observation in (\FPHIone) that, for all $(x,x^*) \in E \times E^*$, $-\infn_{(y,y^*) \in A}\bra{y - x}{y^* - x^*} \le f(x,x^*) - \bra{x}{x^*}$.
\smallskip
(f)\enspace If $x \in \pi_E(\dom\,f)$ then there exists $x^* \in E^*$ such that $f(x,x^*) < \infty$, and so it follows from (e) that there exists $(y,y^*) \in A$ such that $\|y - x\| < 1/n$.   Consequently, $x \in \overline{\pi_E(A)}$.   Thus we have proved that $\pi_E(\dom\,f) \subset \overline{\pi_E(A)}$.   On the other hand, $A \subset \dom\,f$, and so $\overline{\pi_E(A)} = \overline{\pi_E(\dom\,f)}$.   We can prove in an exactly similar way that $\overline{\pi_{E^*}(A)} = \overline{\pi_{E^*}(\dom\,f)}$.   The convexity of the sets $\overline{\pi_E(A)}$ and $\overline{\pi_{E^*}(A)}$ now follows immediately.\qed
\medbreak
In the final results of this section, which are more conveniently stated in terms of multifunctions, we give other consequences of Theorem \EEEQthm.   ``Type (FP)'' (= ``locally maximally monotone'') was defined in \cite\FIVE, Definition 6, p.\ 394\endcite\ and \cite\HBM, Definition 36.5, p.\ 149\endcite, ``type (FPV)'' (= ``maximally monotone locally'') was defined in \cite\FIVE, Definition 7, p.\ 395\endcite\ and \cite\HBM, Definition 36.7, p.\ 150\endcite, ``strongly maximally monotone'' was defined in \cite\FIVE, Definition 8, pp.\ 395--396\endcite\ and  \cite\HBM, Definition 36.9, p.\ 151\endcite, and the statement ``$S + \lambda J_\eta$ is surjective'' was defined in \cite\HBM, (42.2), p.\ 164\endcite.   The facts that strongly representable maximally monotone multifunctions are of type (FP) \big(type (FPV) and strongly maximally monotone, respectively\big) were observed in \cite\VZ, Theorem 22\endcite, \big(\cite\VZ, Remark 6\endcite\ and \cite\VZ, Theorem 23\endcite, respectively\big).   In the above acronyms, ``F'' stands for ``Fiztpatrick'', ``P'' stands for ``Phelps'' and ``V'' stands for ``Veronas''.
\defTheorem \DFPthm
\medbreak
\noindent
{\bf Theorem \DFPthm.}\enspace\slant Let $E$ be a nonzero Banach space, and $S\colon\ E \toto E^*$ be maximally monotone of type (NI). Then:
\smallbreak
\noindent
{\rm(a)}\enspace $S$ is of type (FP).
\smallbreak
\noindent
{\rm(b)}\enspace $S$ is of type (FPV).
\smallbreak
\noindent
{\rm(c)}\enspace $S$ is strongly maximally monotone.   If, further, $S^{-1}\colon\ E^* \toto E$ is coercive, that is to say\break $\inf\bra{S^{-1}x^*}{x^*}/\|x^*\| \to \infty$ as $\|x^*\| \to \infty$, then $D(S) = E$.
\smallbreak
\noindent
{\rm(d)}\enspace For all $\lambda,\eta > 0$, $S + \lambda J_\eta$ is surjective.\endslant
\Proof(a) follows from Theorem \EEEQthm\ and \cite\HBM, Theorem 37.1, pp.\ 153--154\endcite.  (b) follows from Theorem \EEEQthm\ and \cite\HBM, Theorem 39.1, pp.\ 157--158\endcite.   The first assertion in (c) follows from Theorem \EEEQthm\ and \cite\HBM, Theorem 40.1, pp.\ 158--159\endcite, and the second assertion in (c) follows from \cite\HBM, Corollary 41.2, p.\ 160\endcite.   (d) follows from Theorem \EEEQthm\ and \cite\HBM, Theorem 42.8, pp.\ 164\endcite.\qed        
\defSection \FMsec
\bigbreak   
\leftline{\bf \FMsec.\quad Appendix: a nonhausdorff Fenchel--Moreau theorem}
\medskip\noindent
In Theorem \MAXthm, we referred to the Fenchel--Moreau theorem for (possibly nonhausdorff) locally convex spaces.   We shall give a proof of this result in Theorem \FMthm.   When we say that $X$ is a \slant locally convex space\endslant, we mean that $X$ is a nonzero real vector space endowed with a topology compatible with its vector structure and  a base of neighborhoods of $0$ of the form $\big\{x \in X\colon\ S(x) \le 1\big\}_{S \in \Sem(X)}$, where $\Sem(X)$ is a family of seminorms on $X$ such that if $S_1 \in \Sem(X)$ and $S_2 \in \Sem(X)$ then $S_1 \vee S_2 \in \Sem(X)$; and if $S \in \Sem(X)$ and $\lambda \ge 0$ then $\lambda S \in \Sem(X)$.   If $L$ is a linear functional on $X$ then $L$ is continuous if, and only if, there exists $S \in \Sem(X)$ such that $L \le S$ on $X$.
\smallskip
As an example of the construction above, we can suppose that $X$ and $Y$ are vector spaces paired by a bilinear form $\bra\cdot\cdot$.   Then $\big(X,w(X,Y)\big)$ is a locally convex space with determining family of seminorms $\big\{|\bra{\cdot}{y_1}| \vee \cdots \vee |\bra{\cdot}{y_n}|\big\}_{n \ge 1,\ y_1,\dots, y_n \in Y}$.
\smallskip
The author is grateful to Constantin Z\u alinescu for showing him a proof of Theorem \FMthm\ based on the standard (Hausdorff) result and a quotient construction.   The proof we give here is a simplification of the result on \slant Fenchel--Moreau points\endslant\ of \cite\HBL, Theorem 5.3, pp.\ 157--158\endcite\ or \cite\HBM, Theorem 12.2, pp.\ 59--60\endcite\ (which is also valid in the nonhausdorff setting).   
\defTheorem \FMthm
\medbreak
\noindent
{\bf Theorem \FMthm.}\enspace\slant Let $X$ be a locally convex space with defining family of seminorms $\Sem(X)$, and $f \in \PC(X)$ be lower semicontinuous.   Write $X^*$ for the set of continuous linear functionals on $X$.   If $L \in X^*$, define $f^*(L) := \sup_X\big[L - f\big]$.    Let $y \in X$.  Then
$$f(y) = \supn_{L \in X^*}\big[L(y) - f^*(L)\big].\meqno\FMone$$\endslant
\Proof Since, for all $L \in X^*$,
$L(y) - f^*(L) = \infn_{x \in X}\big[L(y) - L(x) + f(x)\big] = (f \episum L)(y)$ and the inequality ``$\ge$'' in (\FMone) is obvious from the definition of $f^*(L)$, we only have to prove that
$$f(y) \le \supn_{L \in X^*}(f \episum L)(y)\big].\meqno\FMtwo$$   
Let $\lambda \in \r$  and $\lambda < f(y)$.   Since $f$ is proper, there exists $z \in \dom\,f$.  Choose $Q \in \Sem(X)$ such that
$$Q(z - x) \le 1 \qlr f(x) > f(z) - 1\meqno\FENMORSUFFone$$
and
$$Q(y - x) \le 1 \qlr f(x) > \lambda.\meqno\FENMORSUFFtwo$$
\par
We first prove that
$$(f \episum Q)(z) \ge f(z) - 1.\meqno\FENMORSUFFthree$$
To this end, let $x$ be an arbitrary element of $X$.  If $Q(z - x) \le 1$ then (\FENMORSUFFone) implies that $f(x) + Q(z - x) \ge f(x) > f(z) - 1$.   If, on the other hand, $Q(z - x) > 1$, let $\gamma := 1/Q(z - x) \in \,]0,1[\,$ and put $u := \gamma x + (1 - \gamma)z$.   Then $Q(z - u) = \gamma Q(z - x) = 1$ and so, from the convexity of $f$, and (\FENMORSUFFone) with $x$ replaced by $u$,
$$\gamma f(x) + (1 - \gamma)f(z) \ge f\big(\gamma x + (1 - \gamma)z\big) = f(u) > f(z) - 1.$$
Substituting in the formula for $\gamma$ and clearing of fractions yields $f(x) + Q(z - x) \ge f(z)$.   This completes the proof of (\FENMORSUFFthree).
\smallskip
Now let $M \ge 1$ and $M \ge \lambda + 2 + Q(z - y) - f(z)$.   We will prove that   
$$(f \episum MQ)(y) \ge \lambda.\meqno\FENMORSUFFfour$$
To this end, let $x$ be an arbitrary element of $X$.  If $Q(y - x) \le 1$ then (\FENMORSUFFtwo) implies that $f(x) + MQ(y - x) \ge f(x) > \lambda$.   If, on the other hand, $Q(y - x) > 1$ then, from (\FENMORSUFFthree),
$$\eqalign{f(x) + MQ(y - x)
&= f(x) + Q(y - x) + (M - 1)Q(y - x)\cr
&\ge f(x) + Q(z - x) - Q(z - y) + (M - 1)\cr
&\ge f(z) - 1 - Q(z - y) + M - 1 \ge \lambda,}$$
which completes the proof of (\FENMORSUFFfour).   The \slant Hahn--Banach--Lagrange theorem\endslant\ of\break \cite\HBL, Theorem 2.9, p.\ 153\endcite\ or \cite\HBM, Theorem 1.11, p.\ 21\endcite\ now provides us with a linear functional $L$ on $X$ such that $L \le MQ$ on $X$ and $(f \episum L)(y) \ge \lambda$, and (\FMtwo) follows by letting $\lambda \to f(y)$.\qed
\bigskip
\leftline{\bf References}
\medskip
\nmbr\BS
\item{[\BS]} R. S. Burachik and B. F. Svaiter, \slant Maximal monotonicity, conjugation and the duality product\endslant,  Proc. Amer. Math. Soc. {\bf 131}  (2003), 2379--2383.
\nmbr\FITZ
\item{[\FITZ]} S. Fitzpatrick, \slant Representing monotone operators
by convex functions\endslant,   Workshop/ Miniconference on
Functional Analysis and Optimization (Canberra, 1988),  59--65, Proc.
Centre Math. Anal. Austral. Nat. Univ., {\bf 20}, Austral. Nat. Univ.,
Canberra, 1988.
\nmbr\GOSSEZ
\item{[\GOSSEZ]}J.- P. Gossez, \slant Op\'erateurs monotones non lin\'eaires dans les espaces de Banach non r\'eflexifs\endslant, J. Math. Anal. Appl. {\bf 34} (1971), 371--395.
\nmbr\GOSSEZC
\item{[\GOSSEZC]}J.- P. Gossez, \slant On a convexity property of the range of a maximal monotone operator\endslant, Proc. Amer. Math. Soc. {\bf 55} (1976), 359--360.
\nmbr\ASTWO
\item{[\ASTWO]}M. Marques Alves and B. F. Svaiter, \slant Br\o ndsted--Rockafellar property and maximality of monotone operators representable by convex functions in non--reflexive Banach spaces.\endslant, J. of Convex Anal., {\bf 15} (2008), 693--706. 
\nmbr\ASTHREE
\item{[\ASTHREE]}-----, \slant A new old class of maximal monotone operators.\endslant, J. of Convex Anal.,  {\bf 16} (2009), 881--890. 
\nmbr\ASD
\item{[\ASD]}-----, \slant On Gossez type (D) maximal monotone operators.\endslant, http://arxiv.org/abs/0903.5332v2, posted April 6, 2009, to appear in J. of Convex Anal., {\bf 17} (2010).
\nmbr\MT
\item{[\MT]}J.--E. Mart\'\i nez-Legaz and M. Th\'era, \slant $\eps$--Subdifferentials in terms of subdifferentials\endslant, Set--Valued
Anal. {\bf 4} (1996), 327--332.
\nmbr\MLT
\item{[\MLT]}-----, \slant
A convex representation of maximal monotone operators\endslant,\
J. Nonlinear Convex Anal. {\bf 2} (2001), 243--247.
\nmbr\PENOT
\item{[\PENOT]}J.--P. Penot, \slant The relevance of convex analysis for the study of monotonicity\endslant,  Nonlinear Anal. {\bf 58}  (2004), 855--871.
\nmbr\PRAGUE
\item{[\PRAGUE]}R. R. Phelps, \slant Lectures on Maximal Monotone
Operators\endslant, Extracta Mathematicae {\bf 12} (1997), 193--230.
\nmbr\FENCHEL
\item{[\FENCHEL]}R. T. Rockafellar, \slant Extension of Fenchel's duality theorem for convex functions\endslant, Duke Math. J. {\bf33} (1966), 81--89.
\nmbr\RTRCA
\item{[\RTRCA]}-----, \slant Convex analysis\endslant, Princeton University Press, Princeton, N.J., 1970.
\nmbr\RANGE
\item{[\RANGE]}S. Simons, \slant The range of a monotone operator\endslant, J. Math. Anal. Appl. {\bf 199} (1996), 176--201.
\nmbr\MANDM
\item{[\MANDM]}-----, \slant Minimax and monotonicity\endslant, 
Lecture Notes in Mathematics {\bf 1693} (1998),\break Springer--Verlag.
\nmbr\BR
\item{[\BR]}-----, \slant Maximal monotone multifunctions of Br{\o}ndsted--Rockafellar type\endslant, Set--Valued Anal. {\bf7} (1999),
255--294.
\nmbr\FIVE
\item{[\FIVE]}-----, \slant Five kinds of maximal monotonicity\endslant, Set--Valued Anal. {\bf9} (2001), 391--409.
\nmbr\PANDM
\item{[\PANDM]}-----, \slant Positive sets and Monotone sets\endslant, J. of Convex Anal., {\bf 14} (2007), 297--317.
\nmbr\HBL
\item{[\HBL]}-----, \slant The Hahn--Banach--Lagrange theorem\endslant, Optimization, {\bf 56} (2007), 149--169.
\nmbr\HBM
\item{[\HBM]}-----, \slant From Hahn--Banach to monotonicity\endslant, 
Lecture Notes in Mathematics, {\bf 1693},\break second edition, (2008), Springer--Verlag.
\nmbr\NRSSD
\item{[\NRSSD]}-----, \slant Nonreflexive Banach SSD spaces\endslant, http://arxiv.org/abs/0810.4579v2, posted\break November 3, 2008.
\nmbr\VZ
\item{[\VZ]}M. D. Voisei and C. Z\u{a}linescu, \slant Strongly--representable operators\endslant, J. of Convex Anal., {\bf 16} (2009), 1011--1033. 
\nmbr\ZAGRODNY
\item{[\ZAGRODNY]}D. Zagrodny, \slant The convexity of the closure of the domain and the range of a maximal monotone multifunction of Type NI\endslant, Set--Valued Anal, {\bf 16} (2008), 759--783.
\nmbr\ZBOOK
\item{[\ZBOOK]} C. Z\u{a}linescu, \slant Convex analysis in general vector spaces\endslant, (2002), World Scientific.
%
%
\Signoff
\bye